\documentclass[12pt,twoside,german]{article}
\usepackage{a4}

\usepackage[ngerman,english]{babel}

\usepackage{amsmath}
\usepackage{amssymb}
\usepackage{amsbsy}
\usepackage{breqn}

\usepackage{caption,longtable,booktabs}
\usepackage{subcaption}
\usepackage[skip=2pt]{caption}

\usepackage{epsfig}
\font\LARGEbsf=cmssdc10 scaled 1700 

\usepackage[pdftex,unicode,colorlinks,linkcolor=blue,citecolor=blue,urlcolor=blue]{hyperref}

\usepackage[boxed,noline,longend]{algorithm2e}
\DontPrintSemicolon
\SetAlCapNameFnt{\small}
\SetAlCapFnt{\small}
\SetAlgoCaptionSeparator{.}
\SetAlCapSkip{1ex}
\SetAlgorithmName{Box}{List of boxes}

\usepackage[labelfont=bf, 
format=plain, 
labelsep=endash, 
justification=raggedright 
]{caption}

\newcommand{\boxe}[1]{\left.#1 \right |_e}
\newcommand{\UPi}{\ten{u}_{\Pi}}
\newcommand{\UPie}{\boxe{\UPi}}
\newcommand{\GradUPi}{\boxe{\nabla\UPi}}{}

\newcommand\blfootnote[1]{%
	\begingroup
	\renewcommand\thefootnote{}\footnote{#1}%
	\addtocounter{footnote}{-1}%
	\endgroup
}

\usepackage{amsmath}
\usepackage{ntheorem}

\usepackage{amssymb}
\usepackage{amsfonts}
\usepackage[]{algorithm2e}
\newcommand {\ten}[1]     {\textbf{\em #1}}

\newcommand {\brm}[1]  {{\bf{#1}}}

\newcommand {\bs}[1] {\mbox{\boldmath{$#1$}}}

\newcommand \rd {d}

\newcommand {\MId}{\ten{1}} 

\newcommand{\tensor}[1]{\left[\begin{matrix}#1\end{matrix}\right]}
\newcommand{\vectorB}[1]{\left(#1\right)}

\newcommand {\PD}[2]   {\frac{\partial #1}{\partial #2}}

\newcommand {\ADGlobalException}[2]   { \left. #1 \right|_{ #2 } }

\newcommand {\IntUd}[3]  {\int\limits_{#1}{#2}{\, \rd #3}}

\def\Grad {\mbox{Grad}}

\def\Div {\mbox{Div}}

\SetKw{KwAnd}{{\bf and}}
\SetKw{KwOr}{{\bf or}}
\SetKw{KwNot}{{\bf not}}
\SetKw{KwFalse}{{\bf false}}
\SetKw{KwTrue}{{\bf true}}
\SetKw{KwPrint}{{\bf print}}
\SetKw{KwStop}{{\bf stop}}

\SetKwComment{KwADBForm}{}{}
\newcommand {\AceGenNamec}[1]    {{\it AceGen#1}}

\theoremstyle{plain}
\theoremstyle{empty}
\newtheorem{plaincode}{}
	\ifx \copytomathematica \undefined
	
	\fi
	
\newcommand\exceq{\stackrel{\smash{\scriptscriptstyle\mathrm{!}}}{=}}
\font\LARGEbsf=cmssdc10 scaled 2100      

\begin{document}
%
\begin{figure}[t]
  \includegraphics[width=1.0\textwidth]{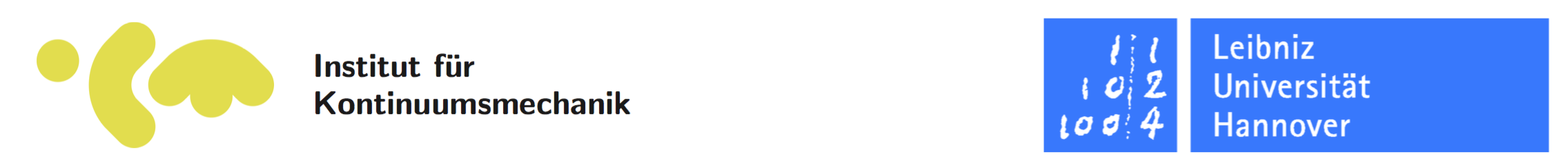}
\end{figure}

\phantom{ii}
\vspace{40 mm}
\begin{center}
{\Large \bf Virtual Element Formulation \\[2mm] 
For Finite Strain Elastodynamics}\\ \vspace{10 mm}
Mertcan Cihan, Fadi Aldakheel, Bla\v{z} Hudobivnik, Peter Wriggers \\ \vspace{10 mm}
{\bf Preprint}\\
{\it February 06, 2020}
\end{center}

\newpage

\begin{center}
{\Large \bf Virtual Element Formulation \\[2mm] 
For Finite Strain Elastodynamics}\\ \vspace{10 mm}
\end{center}

\begin{center}
Mertcan Cihan${}^{\ast}$, Fadi Aldakheel, Bla\v{z} Hudobivnik, Peter Wriggers 

Institute for Continuum Mechanics, Leibniz Universit\"at Hannover, Germany
\end{center}

\blfootnote{$^*$ Corresponding author.\\
\hspace*{9mm}\textit{E-mail addresses:} cihan@ikm.uni-hannover.de (Mertcan Cihan), aldakheel@ikm.uni-hannover.de (Fadi Aldakheel), hudobivnik@ikm.uni-hannover.de (Bla\v{z} Hudobivnik), wriggers@ikm.uni-hannover.de (Peter Wriggers).
\vspace{6mm}\\
\hspace{-4mm} \textit{Preprint}}

\vspace{4mm}
\begin{center}
{\bf \large Abstract}
\bigskip
\vspace{1cm}
{\footnotesize
\begin{minipage}{14.5cm}
\noindent

This work provides an efficient virtual element scheme for the modeling of {\it nonlinear} elastodynamics undergoing large deformations. The virtual element method (VEM) has been applied to various engineering problems such as elasto-plasticity, multiphysics, damage and fracture mechanics. This work focuses on the extension of VEM towards {\it dynamic applications}. Within this framework, we employ low-order ansatz functions in one, two and three dimensions that having arbitrary convex or concave polygonal elements. The formulations considered in this contribution are based on minimization of potential function for both the static and the dynamic behavior. While the stiffness-matrix needs a suitable stabilization, the mass-matrix can be calculated using only the {\it projection part}. For the implicit time integration scheme, Newmark-Method is used. To show the performance of the method, various numerical examples in 1D, 2D and 3D are presented.  
\end{minipage}
}
\end{center}

\section{Introduction}
The virtual element method (VEM) can be seen as an extension of the classical finite element method (FEM) based on Galerkin projection. It allows meshes with highly irregular shaped elements, including non-convex shapes, as outlined in \cite{veiga+etal13,wriggers19}. This gives more flexibility and new possibilities to geometry discretization in solid- and fluid-mechanics. The large number of positive properties of VEM increases the variety of possible applications in engineering and science. Recent works on virtual elements have been devoted to linear elastic deformations in \cite{veiga+brezzi+etal13,gain+etal14,artioli17}, contact problems in \cite{wriggers+rust+reddy16}, elasto-plastic deformations in \cite{hudobivnik2019low,aldakheel-themero19,ArVeLoSa17b}, anisotropic materials in  \cite{wriggers+hudobivnik+korelc18,wriggers+hudobivnik+schroeder18,Reddy2019}, curvilinear virtual elements for 2D solid mechanics applications in \cite{artioli2019curvilinear}, hyperelastic materials at finite deformations in \cite{ChVePa16,wriggers+reddy+rust+hudobivnik17}, crack-propagation for 2D elastic solids at small strains in \cite{Hussein+aldakheel18} and phase-field modeling of brittle and ductile fracture in \cite{aldakheel+blaz+wriggers18,aldakheeletal18}.
\\
Despite the fact that dynamic behavior has a strong influence on the mechanical properties and the prediction of their real response, most of the investigations introduced above are only done for static problems so far. Thus the element mass-matrix is needed to be calculated. In this regard, \cite{Park19} proposed a virtual element method for {\it linear} elastodynamics problems. However their formulations are restricted to small strain setting, hence it is not appropriate for large deformations. This has motivated the presented contribution to extend the application of VEM from the static to the dynamic case in the finite deformation range.\\
Typically the construction of a virtual element is divided into a projection step and a stabilization step. Within the projection step, a quantity $\varphi_h$ is replaced by its {projection $\varphi_\Pi$} onto a polynomial space. Using this projected quantity in the weak formulation or energy functional yields to a rank-deficient structure which needs to be stabilized. In the second step, the stabilization term, which is a function of the difference $\varphi_h - \varphi_\Pi$ between the original variable and the projected quantity needs to be evaluated. There are various possibilities to evaluate this stabilization term. To this end, Da Veiga et al. \cite{BeLoMo15} proposed a stabilization term, where all integrations take place on the element boundaries. Wriggers et al. presented in \cite{wriggers+rust+reddy16} a novel stabilization technique, which was first described for finite elements in Nadler and Rubin \cite{Nadl03}, generalized in Boerner et al. \cite{Boer07} and simplified in Krysl \cite{Krys16} for the stabilization procedure of a mean-strain hexahedron. In this framework, the integration is carried out over a triangulated sub-mesh, which uses the same nodes as the original mesh. The method presented in this contribution based on the stabilization technique of \cite{wriggers+rust+reddy16}. In order to model the dynamic behavior of the body, we define a specific potential function, where the second derivative of it with respect to the global unknowns yields the mass-matrix. As a key advantage of this approach in comparison to \cite{Park19}, only the stiffness-matrix needs to be stabilized, whereas the mass-matrix is only computed using the {\it projection part}. Hence no stabilization is required for the mass-matrix as in the case of \cite{Park19}. For the time integration scheme, we utilized the implicit Newmark method as documented in \cite{Newm59, Wood80}.
\\
The structure of the presented work is as follows. In Section \ref{elasticity} the governing equations for nonlinear elastodynamics are outlined. Section \ref{vem_base} summarizes the virtual element formulation. It includes details on the computation of the element mass-matrix. To verify the proposed virtual element formulations, a various number of examples are demonstrated and discussed in Section \ref{Examples}. Section \ref{Summary} briefly summarizes the work and gives some concluding remarks.

\section{Governing equations for Nonlinear Elastodynamics}
\label{elasticity}
\begin{figure}[!tbh]
    \centering
    \includegraphics[width=0.45\textwidth]{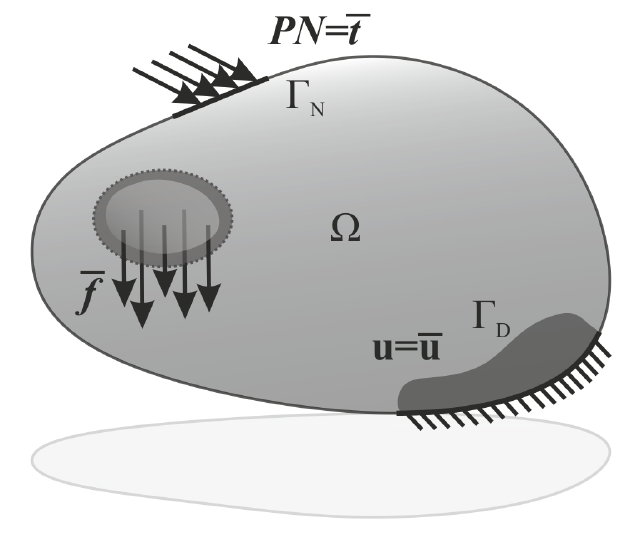}
\caption{\label{fig:continuum} Solid with boundary conditions}
\end{figure}

In this section we summarize the finite strain elasto-static formulation (see e.g. \cite{Wri:2008:nfe,KoWr16,KoSt14}) and supplement it by the dynamic effect. For that consider an elastic Body $\Omega\subset \mathbb{R}^3$ with boundary $\Gamma$. This boundary is decomposed into a non-overlapping Dirichlet $\Gamma_{D}$ and Neumann $\Gamma_{N}$ boundary conditions such that $\Gamma_{D} \cup \Gamma_{N} = \Gamma$, see Figure \ref{fig:continuum}.

The position $\ten{x}$ of a material point in the current configuration is given by the deformation map
\begin{align}
\ten{x} = \bs{\varphi}(\ten{X},t) = \ten{X} + \ten{u}(\ten{X},t)\,,
\label{motion}
\end{align}
where $\ten{X}$ is the position of a material point in the initial configuration and $\ten{u}(\ten{X},t)$ is the displacement. In the further course of this work we will skip the explicit specification of the dependence of variables on the initial configuration and time thus we will write: $\ten{u}=\ten{u}(\ten{X},t)$. In order to transform quantities which are defined with respect to the deformed configuration to the reference configuration and vice versa, we define the deformation gradient $\ten{F}$ as
\begin{equation}
\ten{F} = \Grad\,\bs{\varphi} = \nabla_{\ten{X}}\, \bs{\varphi}\,,
\label{F}
\end{equation}
where the gradient is evaluated with respect to the initial configuration $\ten{X}$. We further define the right Cauchy-Green tensor $ \ten{C}(\ten{u} )$ with $ \ten{C}= { \ten{F}}^T \ten{F}$ as a strain measure and the Jacobian $J(\ten{u} )$ with $J= \det  \ten{F}$ as a volume map.\\

The solid $\Omega$ has to satisfy the balance of linear momentum 
\begin{equation}
\Div\,\ten{P} + \overline{\ten{f}} = \rho \ddot{\ten{u}} 
\quad \mbox{with} \quad
\ten{P} = \ten{F\,S} \; ,
\label{equil:balance of linear momentum}
\end{equation}
where $\overline{\ten{f}}$ are the body forces and $\ten{P}$, $\ten{S}$ are the $1^{st}$ and $2^{nd}$ Piola-Kirchhoff stresses, respectively. The right side of the equation \eqref{equil:balance of linear momentum}$_1$ is taking the dynamic effects $\rho \ddot{\ten{u}}$ into consideration. The Dirichlet and Neumann boundary conditions are defined by
\begin{align}
\ten{u} & = \bar{\ten{u}} \qquad \mbox{on}\ {\Gamma_{D}}\,,  \label{Dbc} \\
\ten{P}\ten{N} & = \bar{\ten{t}}\qquad \mbox{on}\ {\Gamma_{N}}\,, \label{Nbc} 
\end{align}
here $\ten{N}$ is the outward unit normal vector related to the initial configuration, $\bar{\ten{u}}$ represents the prescribed displacement on the Dirichlet boundary $\Gamma_{D}$, and $\bar{\ten{t}}$ depicts the surface traction at the Neumann boundary  $ \Gamma_{N}$, as illustrated in Figure \ref{fig:continuum}. The weak formulation of the elastodynamics problem in \eqref{equil:balance of linear momentum} then takes the form
\begin{equation}
G(\ten{u},\delta \ten {u}) = \int\limits_{\Omega} \left  [\ten{S}(\ten{u}) : \frac{1}{2}\ten{C}(\delta \ten u) - \overline{\ten{f}}\cdot \delta \ten u\ + \rho \ddot{\ten{u}} \cdot \delta \ten u \right ] \ d\Omega -   \int\limits_{\Gamma_N} \bar{\ten{t}}\cdot \delta \ten u \ d\Gamma\, ,
\label{WeakF}
\end{equation}
where $\delta \ten u$ is the test function of the displacement $\ten u$. A homogeneous compressible isotropic elastic material is considered, here we use the Neo-Hookean strain energy function 
\begin{equation}
\Psi= \frac{\kappa}{4}(I_3-1- \ln {I_3})  + \frac{\mu}{2} (I_3^{-1/3}I_1 -3) \; ,
\label{Diso}
\end{equation}
in terms of the bulk $\kappa$ and shear $\mu$ modulus and the invariants of the right Cauchy-Green tensor $I_1=\mbox{tr} \,\ten{C} $, $I_3=\det \ten{C}$.\\ 
\\
\newline
With the above set of equations, the finite strain elastodynamic problem is well formulated. Next, we use the potential function as a starting point for the development of a discretization method. The static part of the potential is defined as
\begin{equation}
U^{stat}(\ten{u}) = \int\limits_{\Omega} \left  [\Psi(\ten{u}) -   \overline{\ten{f}}\cdot \ten{u}\ \right ] \ d\Omega -   \int\limits_{\Gamma_N} \bar{\ten{t}}\cdot \ten{u} \ d\Gamma\, ,
\label{PsPotEnergStat}
\end{equation}
whereas the dynamic part of the potential is the kinetic energy that describes inertial effects takes the form 
\begin{equation}
\mathcal{K}(\ten{u}) =\frac{1}{2} \int\limits_{\Omega} \rho  {\dot{\ten{u}}}^2 \ d\Omega\,,
\label{PsPotEnergDyn}
\end{equation}
where $\rho$ is the density of the solid.

\section{Formulation of the Virtual Element Method}
\label{vem_base}

The main idea of the virtual element method is to use a Galerkin projection, which maps the primary fields (displacement in this work) to a specific polynomial ansatz space. Thus, the domain $\Omega$ can be discretized with non-overlapping polygon in 2D or polyhedral elements in 3D which do not need to have convex shapes \cite{veiga+brezzi+etal13}. Since VEM has no isoparametric mapping, the ansatz functions are given in terms of the coordinates $\ten{X}$ in the initial configuration. Here the ansatz functions for the virtual element are based on linear functions, therefore the nodes are placed at the element vertices.

\subsection{VEM Ansatz}
\label{decomp}
In general, for finite strains the deformation map $\bs\varphi = \ten{ X} + \ten{u}$ has to be discretized. But as the coordinates $\ten{X}$ in the initial configuration are exactly known, we can reduce the discretization to the displacement field $\ten u= u_i\,\ten{E}_i$ where $\ten{E}_i$ are the basis vectors with respect to the initial configuration in the three-dimensional space $i\in \{1,2,3\}$. 

The central concept of the virtual element method relies on the split of the ansatz space $ \ten{u}_h $ into a projected part $\ten{u}_{\Pi} $ and a remainder $ \ten{u}_h - \ten{u}_{\Pi}$ as
\begin{equation}
 \ten{u}_h =  \ten{u}_{\Pi}+ ( \ten{u}_h - \ten{u}_{\Pi})
 \label{eq:usplit}
\end{equation}
For a linear ansatz, the projection $\ten{u}_{\,\Pi}$ at element level takes for three-dimensional elements the form
\begin{equation}
\begin{array}{r@{\ }c@{\ }l}
\ten{u}_{\,\Pi} &=& \left(\brm{N}_{\,\Pi}\cdot \brm{a}_i\right)\,\ten{E}_i \, ,\\[2mm] 
\brm{N}_{\,\Pi} & =& \vectorB{1, X, Y, Z}\, ,\\[2mm]
\brm{a}_i  & =& \vectorB{a_{i\,1}, a_{i\,2}, a_{i\,3}, a_{i\,4}}\, ,
\end{array}
\label{eq:projupi3D}
\end{equation}
where $\brm{a}$ represents the twelve unknown virtual parameter $\brm{a}=\bigcup\brm a_i$ which have to be determined. Instead of using the polynomial $\brm{N}_{\Pi} $ in equation \eqref{eq:projupi3D} as the interpolation function, a scaled ansatz can be used, for details see e.g. \cite{artioli17}.
Furthermore, the projection $\ten{u}_{\Pi}$ has to fulfill the orthogonality condition, as discussed in the work of \cite{beirao2013c}. Hence $\nabla \ten{u}_{\Pi}$ is computed through the Galerkin projection as
\begin{equation}
\int_{\Omega_{e}} \nabla \ten{p} \cdot (\nabla \ten{u}_{\Pi} -\nabla \ten{u}_h) \ d\Omega = 0 \, ,
\label{eq:Galerkin Projection}
\end{equation}
where $\ten{p}$ is a polynomial function which has been chosen similarly to the projection $\ten{u}_{\,\Pi}$, see \eqref{eq:projupi3D}. Since we use linear ansatz functions, $\nabla \ten{p}$ and $\nabla \ten{u}_{\Pi}$ are constant and can be shifted out of the integral as
\begin{equation}
\nabla \ten{u}_{\Pi} = \frac{1}{\Omega_e} \int_{\Omega_{e}} \nabla  \ten{u}_h \ d\Omega\,.  
\label{eq:orthogonality_condition2}
\end{equation}
Applying integration by parts to \eqref{eq:orthogonality_condition2}, we obtain
\begin{equation}
\GradUPi \exceq \frac{1}{\Omega_e}\int_{\Gamma_{e}} \ten{u}_h \otimes \ten{N}\,\ d \Gamma \ ,
\label{eq:cProj}
\end{equation}
at the element level. Here $\ten{N}$ denotes the normal vector on the reference boundary $\Gamma_e$ of the domain $\Omega_e$, which belongs to a virtual element $e$. Element quantities, which have constant values within the entire element $e$, are denoted by $\boxe{\square}$.
With this simplification the projection $\ten{u}_{\Pi}$ is defined as outlined in \cite{wriggers+reddy+rust+hudobivnik17}.

\begin{figure}[!h]
\centering
\scalebox{1}{\includegraphics[width=0.9\textwidth]{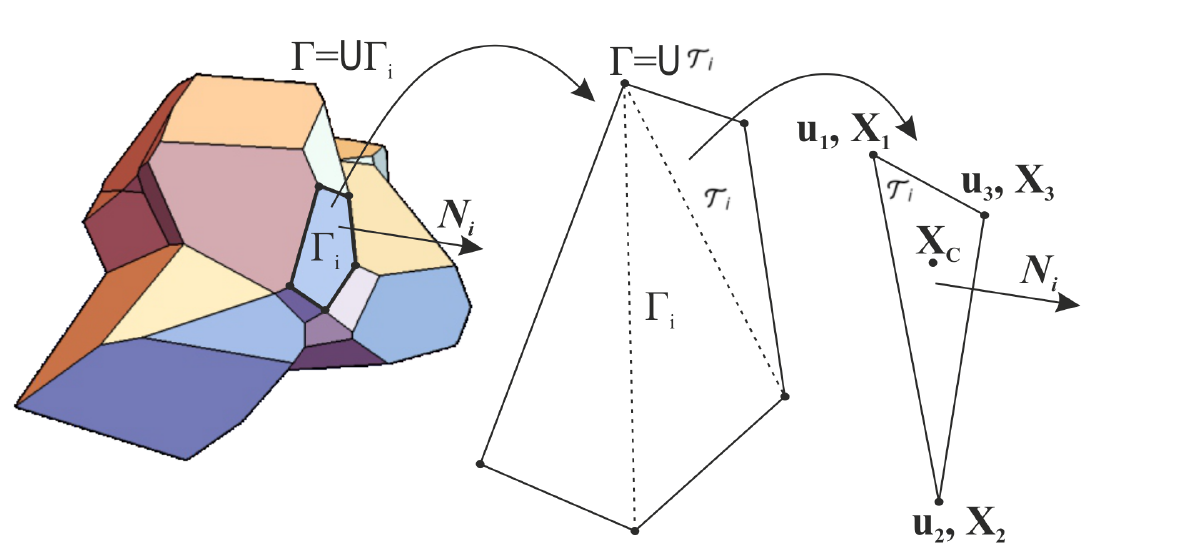}}
\caption{\label{VEM_EL_Face}Virtual element faces split into multiple triangles}
\end{figure}

By employing the linear ansatz space, the left hand side of (\ref{eq:cProj}) takes the simple form
\begin{equation}
\GradUPi = \tensor{
a_{1\,2} & a_{1\,3} & a_{1\,4} \\
a_{2\,2} & a_{2\,3} & a_{2\,4} \\
a_{3\,2} & a_{3\,3} & a_{3\,4}} \,.
 \label{eq:cGrad}
\end{equation}
In the 2D case, the right hand side of \eqref{eq:cProj} is evaluated along the edges. As the displacements are known at the boundary, which are straight line segments, a linear ansatz for the displacements is used, see \cite{wriggers+reddy+rust+hudobivnik17}. However, in the 3D case, the element boundary consists of polygonal faces. Therefore the evaluation of the integral in \eqref{eq:cProj} is not straight forward, unless an appropriate ansatz is found. For the evaluation, there are two possible methods available. The first one is presented in \cite{gain+etal14}, where the faces are subdivided in quadrilateral elements where the corners of the quadrilateral elements have certain positions. Finally the evaluation of the integral is carried out on those quadrilateral elements.
An alternative option is to subdivide the element faces into 3 noded triangles. The integration is then carried out over the triangles of the polygonal faces by using the standard ansatz function for a linear triangle and Gauss integration:
\begin{align}
\brm N_h & = \vectorB{\xi, \eta, 1 - \xi - \eta}  \label{eq:NhP1} \\
\ten u_h^{\mathcal{T}} & = \brm N_h  \brm u_I \qquad \forall I\in \mathcal{T} \, , \label{eq:uhP1}
\end{align} 
as outlined in \cite{wriggers19}. Here $\ten u_h^{\mathcal{T}} $ denotes the linear ansatz for the displacements at each triangle of the polygonal faces.
$\brm u_I$ is a list which contains the three nodal displacement vectors of the triangle $\mathcal{T}$. $\xi$ and $\eta$ are the local dimensionless coordinates at the element level. The local nodes of $\mathcal{T}$ and the outward normal vector $\emph{\textbf{N}}_i$ are visualized in Figure \ref{VEM_EL_Face}.
Finally we are able to compute the right hand side of \eqref{eq:cProj}. Using \eqref{eq:uhP1}, the integral in \eqref{eq:cProj} takes the form:
\begin{equation}
\frac{1}{\Omega_e}\int\limits_{\Gamma_{e}}{\ten{u}_h \otimes \ten{N}}\,d{\Gamma} 
= \frac{1}{\Omega_e}  \sum_{k=1}^{n_f}\int\limits_{\Gamma_{k}}{\ten{u}_h^{\mathcal{T}}  \otimes \ten{N}_k}\,d {\Gamma} 
= \frac{1}{\Omega_e}  \sum_{k=1}^{n_f} \sum_{g=1}^{n_g} w_g N_\zeta \ten{u}_{h\,g}^{\mathcal{T}}  \otimes \ten{N}_g
 \label{eq:cGradh}
\end{equation}
Here $n_f$ is the number of element faces. For an integration over triangles with linear shape functions \eqref{eq:NhP1} one point quadrature with $n_g=1$ Gauss point and $w_g = 1/2$ Gauss weight is sufficient. $N_\zeta$ is the Jacobian of transformation from the  reference to the initial configuration. 
$\square_g$ denotes quantities which are evaluated at the Gauss point with the local coordinates $\xi=1/3$ and $\eta=1/3$. The normal vector $\ten{N}$ and the  Jacobian of the isoparametric mapping $N_\zeta$ are evaluated as follows:
\begin{align}
\ten X^{\mathcal{T}} & = \brm N_h \,\brm X_I \qquad {\forall I\in \mathcal{T}},\label{eq:X}\\
\ten g_\xi &= \PD{\ten X^{\mathcal{T}}}{\xi},\,\, \ten g_\eta= \PD{\ten X^{\mathcal{T}}}{\eta},\,\, \ten g_\zeta= \ten g_\xi \times \ten g_\eta\,,\\
N_\zeta &=|\ten g_\zeta|, \, \ten N= \frac{\ten g_\zeta}{N_\zeta}\,.
\end{align}
All quantities are related to the initial configuration. 
\\
\\
Comparing \eqref{eq:cGrad} and \eqref{eq:cGradh} we obtain the unknown virtual parameters $\ADGlobalException{\brm a}{i\in(4,...,12)}$ by inspection, for further details see e.g. \cite{wriggers+reddy+rust+hudobivnik17}. Since only the gradient of the projection $\GradUPi$ is needed to define the strain energy function of the static part, we actually do not have to compute the virtual parameters which are related to the constant parts. However, we will later see that the projected displacements $\ten{u}_{\Pi} $, including the virtual parameters which are related to the constant parts, are needed to construct the mass-matrix. Thus we need to supplement our formulation by a further condition to obtain the constants $\ADGlobalException{\brm a}{i\in(1,2,3)}$. For this purpose we adopt the condition (see for example \cite{veiga+brezzi+etal13}) that the sum of the nodal values of $\ten{u}_h$ and of its projection $\ten{u}_{\Pi} $ are equal over the entire domain. This yields for each virtual element $\Omega_{e}$ to the following condition
\begin{equation}
\frac{1}{n_{V}}\sum^{n_{V}}_{I=1} \ten{u}_{\Pi} (\brm{X}_I) =\frac{1}{n_{V}} \sum^{n_{V}}_{I=1}\ten{u}_{h}(\brm{X}_I) \,, \label{eq:avefn} 
\end{equation}
where $n_V$ is the number of boundary nodes and $\brm{X}_I$ the initial coordinates of the nodal point $I$. The sum includes all boundary nodes $n_V$. 
By substituting \eqref{eq:projupi3D} in \eqref{eq:avefn} we can express the missing three parameters in terms of the nodal displacements and the already known virtual parameters $\ADGlobalException{\brm a}{i\in(4,...,12)}$:
\begin{equation}
{\vectorB{a_{1\,1}, a_{2\,1}, a_{3\,1}}} = \frac{1}{n_V} \sum^{n_{V}}_{I=1} \left(\brm{u}_I - \GradUPi \, \brm X_I \right)\label{eq:translation}
\; .
\end{equation}
Finally with equation \eqref{eq:cProj} and \eqref{eq:translation} the ansatz function $\ten{u}_{\Pi} $ of the virtual element is completely defined in terms of the element nodal displacements $\brm{u}_e$:
\begin{equation}
\brm{u}_e= \left(\brm u_I\right)_{\forall I\in \lbrace 1,..., n_V\rbrace} = \lbrace \brm u_1, \brm u_2, ... , \brm u_{n_V} \rbrace.
\label{eq:nodalU}
\end{equation}

\subsection{Implicit time integration}
For the time integration scheme, we use the implicit Newmark method outlined in \cite{Newm59, Wood80}. The equations for the velocity $\dot {\brm{u}}=\ten{v}$ and acceleration $\ddot {\brm{u}}$ at time step $t_{n+1}$ are given as

\begin{equation}
\dot {\brm{u}}_{n+1}( {\brm{u}}_{n+1})=\dfrac{\gamma}{\zeta \Delta t} \left( {\brm{u}}_{n+1}-{\brm{u}}_n \right)-\left(\dfrac{\gamma}{\zeta}-1 \right) \dot {\brm{u}}_n- \left(\dfrac{\gamma}{\zeta}-1 \right) \Delta t \ddot {\brm{u}}_n
\label{eq:velocity-Newmark}
\end{equation}
\begin{equation}
\dot {\brm{u}}_{n+1}( {\brm{u}}_{n+1})= \dfrac{1}{\zeta \Delta t^2} ({\brm{u}}_{n+1}-{\brm{u}}_n)-\dfrac{1}{\zeta \Delta t} \dot {\brm{u}}_n - \left( \dfrac{1}{2 \zeta}-1 \right) \ddot {\brm{u}}_n \, ,
\label{eq:acceleration-Newmark}
\end{equation}
where the Newmark parameters are chosen as $\zeta=1/4$ and $\gamma=1/2$.


\subsection{Construction of the element mass-matrix for VEM}
\label{Mass Matrix}
For the calculation of the mass-matrix, we start from the dynamic potential in equation \eqref{PsPotEnergDyn}. From \eqref{eq:cProj} and \eqref{eq:nodalU} we obtain the virtual parameters $\brm{a}$ in terms of the unknown displacements $\brm{u}_e$ as 
\begin{equation}
\ten{a} = \brm{\tilde{\Pi}}^\nabla  \brm{u}_e \ ,
\label{eq:VEM Ansatz pii}
\end{equation}
where $\brm{\tilde{\Pi}}^\nabla$ is a constant matrix. Hence the projection $\ten{u}_{\,\Pi}$ can be rewritten as
\begin{equation}
\ten{u}_{\,\Pi} = \brm{H}  \, \ten{a} =\brm{H}\, \brm{\tilde{\Pi}}^\nabla  \brm{u}_e  \, ,
\label{eq:coupling virtual parameters global unknowns}
\end{equation}
where $\brm{H}$ is the matrix representation of the ansatz functions $\brm{N}_{\,\Pi}$. It is defined in the three-dimensional case as

\begin{equation}
\begin{array}{c@{\ }c@{\ }c}
\brm{H}&=& 
\begin{bmatrix}
1&0&0&X&0&0&Y&0&0&Z&0&0\\
0&1&0&0&X&0&0&Y&0&0&Z&0\\
0&0&1&0&0&X&0&0&Y&0&0&Z\\
\end{bmatrix}\\
\end{array} \, .
\label{eq:H Matrix}
\end{equation}

The matrix representation of the projection operator $\brm{\tilde{\Pi}}^\nabla$ in \eqref{eq:coupling virtual parameters global unknowns} 
can be derived with the help of \eqref{eq:cGradh} and \eqref{eq:translation}, as discussed in \cite{wriggers+reddy+rust+hudobivnik17}.

To compute a good approximation for the mass-matrix, we use the projection $\ten{u}_{\Pi}$ which is now defined by the orthogonality condition:
\begin{equation}
\int_{\Omega_{e}} \ten{p} \cdot \ten{u}_{\Pi}  \ d\Omega = \int_{\Omega_{e}} \ten{p} \cdot \ten{u}_h  \ d\Omega\, .  
\label{eq:Orthogonality condition L2}
\end{equation}
Usually one has to compute a new projection which is not the same as for the construction of the stiffness matrix. Nevertheless, it has been shown in \cite{beirao2013c, Ahma13} that it is possible to use the same projection operator $\brm{\tilde{\Pi}}^\nabla$ for $k \leq 2$ where $k$ is the order of the ansatz functions and it does not lead to an error. Thus we use the same projection operator $\brm{\tilde{\Pi}}^\nabla$, which has been used for the static part as

\begin{equation}
\ddot {\ten{u}}_{\Pi} =\brm{H} \brm{\tilde{\Pi}}^\nabla  \ddot{\brm{u}}_e \, .
\label{eq:VEM Ansatz ddotu}
\end{equation}
For the unknown accelerations $\ddot{\ten{u}}_h$ we are using the same split as for the unknown displacements in equation \eqref{eq:usplit}, yields
\begin{equation}
\ddot {\ten{u}}_h = \ddot {\ten{u}}_{\Pi} + (\ddot {\ten{u}}_h - \ddot {\ten{u}}_{\Pi})\, .
\label{eq:asplit}
\end{equation}
For the construction of the elastodynamic virtual element, it is computationally advantageous to employ the software tool \AceGenNamec, see \cite{KoWr16}. It provides the most efficient element routines when a potential formulation is used. Thus we construct a {\it specific} pseudo-potential for inertia term
\begin{equation}
U^{dyn}(\ten{u}) = \int\limits_{\Omega} \rho  \ddot{\ten{u}}\, \cdot \ten{u} \ d\Omega\, ,
\label{new-dyn}
\end{equation}
where the first variation has to be performed for fixed $\ddot{\ten{u}}$.
The total potential function is now defined as
\begin{equation}
U=U^{stat}+U^{dyn} \, ,
\label{eq:full potential}
\end{equation}
which depends on both the static and inertia parts. The variation of \eqref{eq:full potential} yields exactly the weak form of \eqref{equil:balance of linear momentum} when considering the nonlinear dependency of the $2^{nd}$ Piola-Kirchhoff stress $\ten{S}$ on the displacement $\ten{u}$. Therefore the usage of the pseudo potential is equivalent to using the weak form \eqref{WeakF} directly.
\\
\newline
If we insert both equations \eqref{eq:usplit} and \eqref{eq:asplit} for the displacements and the accelerations in \eqref{new-dyn}, we obtain:
\begin{equation}
U^{dyn}(\ten{u}) = \int\limits_{\Omega} \rho \ddot {\ten{u}}_\Pi \cdot \ten{u}_\Pi \ d\Omega + \int\limits_{\Omega} \rho  (\ddot {\ten{u}}_h - \ddot {\ten{u}}_\Pi) \cdot (\ten{u}_h -\ten{u}_\Pi) \ d\Omega \, ,
\label{eq:Dynamic_Potential}
\end{equation}
where coupled terms vanish due to the orthogonality condition \eqref{eq:Orthogonality condition L2}. The first term in \eqref{eq:Dynamic_Potential} is the consistency part, whereas the second term is the stabilization part. It is sufficient to use the consistency term for the construction of the mass-matrix, without any stabilization, when the problem is without any reaction term, as shown in \cite{beirao2013c, Ahma13}. To compute the mass-matrix from $U^{dyn}$, we need to compute the first and second derivative with respect to the global unknowns. By utilizing the relationship between the projected values and the unknown values for the displacement and the accelerations, the following expression for the modified dynamic potential function results
\begin{equation}
U^{dyn}(\ten{u}_\Pi) = \int\limits_{\Omega} \rho \ddot {\ten{u}}_\Pi \cdot \ten{u}_\Pi \ d\Omega = \int\limits_{\Omega} \rho \ \brm{a}^T \brm{H}^T \brm{H} \ddot{\brm{a}} \ d\Omega
\label{eq:Dynamic_Potential2}
\end{equation}
Hereby the virtual accelerations and displacements are constant over the entire domain of the element therefore they can be shifted out of the integral. Thus the integral yields
\begin{equation}
U^{dyn}(\ten{u}_\Pi) = \brm{a}^T \Big[\int\limits_{\Omega} \rho \brm{H}^T \brm{H} \ d\Omega \Big] \ddot {\brm{a}} \quad \mbox{with } 
\begin{array}{c@{\ }c@{\ }c}
\brm{H}^T \brm{H}&=& 
\begin{bmatrix}
1&0&X&0&Y&0\\
0&1&0&X&0&Y\\
X&0&X^2&0&XY&0\\
0&X&0&X^2&0&XY\\
Y&0&XY&0&Y^2&0\\
0&Y&0&XY&0&Y^2\\
\end{bmatrix}\\
\end{array}
\label{eq:Integral_HTH}
\end{equation}
This integral can be evaluated in different ways which will be explained in Section \ref{Consistency part}.
\\
\newline
The first derivative of the dynamic potential $U^{dyn}$ is computed holding the acceleration $\ddot {\ten{u}}_e$ constant:
\begin{equation}
\ten{R}^{dyn} = \left. \dfrac{\partial U^{dyn}(\UPie)}{ \partial \brm{u}_e} \right|_{\ddot {\brm{u}}_e = const.}=\brm{M} \cdot \ddot {\brm{u}}_e
\quad \mbox{with} \quad
\brm{M}=\left(\brm{\tilde{\Pi}}^\nabla\right)^T \int\limits_{\Omega} \rho \brm{H}^T \brm{H} \ d\Omega \ \brm{\tilde{\Pi}}^\nabla
\label{eq:1st variation}
\end{equation}
Before computing the second derivative, the Newmark method is used for the implicit time integration. With \eqref{eq:acceleration-Newmark}, the residual for the dynamic part follows as:
\begin{equation}
\begin{aligned}
\ten{R}^{dyn} = \brm{M} \cdot \left[\dfrac{1}{\zeta \Delta t^2} (\brm{u}_{e,n+1}-\brm{u}_{e,n})-\dfrac{1}{\zeta \Delta t} \dot {\brm{u}}_{e,n} -(\dfrac{1}{2 \zeta}-1)\ddot {\brm{u}}_{e,n} \right]
\end{aligned}
\label{eq:1st variation with NM}
\end{equation}
The second derivative of $U^{dyn}$ leads then to the dynamic part of the tangent 
\begin{equation}
\dfrac{\partial^2 U^{dyn}(\UPie)}{ \partial \brm{u}_e^2} = \brm{M} \cdot \dfrac{1}{\zeta \Delta t^2} 
\label{eq:2nd variation}
\end{equation}

\subsection{Construction of the virtual element}
As introduced in Section \eqref{vem_base}, the formulation of a virtual element undergoing large deformations is based on a split of the energy into a constant part and an associated stabilization term. The nodal degrees of freedom of an element are approximated with one interpolation function per coordinate direction in each element. Thus the consistency part does not lead to a stable formulation and a stabilization term is required. The idea of stabilizing the formulation is analogous to the stabilization of the classical finite elements with reduced integration, developed by \cite{Krys16}.
For the construction of the virtual element method we start from the potential function \eqref{eq:full potential}. After summing up all element contributions for the $\ten{n}_e$ virtual elements we obtain the following expression:

\begin{equation}
U(\ten{u}) = \raise3pt
\hbox{$
	\hbox{\scriptsize $n_e$}\atop{\hbox{\LARGEbsf A}\atop
		{\scriptstyle e=1}}$} \left [ U_c(\UPie) + U_{stab} (\boxe{\ten{u}_h} -\UPie)\,\right]
\label{eq:Pseudopotential stat+dyn proj+stab}
\end{equation}

\subsubsection{Consistency part}
\label{Consistency part}
For the consistency part we use the projection $\ten{u}_\Pi$ as introduced in Section \ref{decomp} and therefore the first part of equation \eqref{eq:Pseudopotential stat+dyn proj+stab} for each element is given by
\begin{equation}
U_c(\left . \UPi \right |_e) = \IntUd{\Omega_e}{\left[\Psi(\UPi) - \overline{\ten{f}} \cdot \UPi \right]}{\Omega}  -   \IntUd{{\Gamma}_e^N}{\bar{\ten{t}}\cdot \UPi}{\Gamma} + \int\limits_{\Omega_e} \rho  \ddot{\ten{u}}_\Pi \cdot \ten{u}_\Pi \ d\Omega
\label{PotEnerg_c}
\end{equation}
The gradient of the projection $ \left.\nabla \ten{u}_{\Pi} \right |_e$ is constant on the entire domain $\Omega_e$ thus all kinematic quantities such as $ \left. \ten{F}\right |_e=  \MId + \GradUPi$ are constant as well. Hence the integration of the strain energy function can be simplified as:
\begin{equation}
\IntUd{\Omega_e}{\Psi( \left. \ten{C} \right |_{e})}{\Omega} =
\Psi( \left. \ten{C} \right |_{e})  \,\Omega_e \ ,
\label{eq:vem_stat_const}
\end{equation}
which is still nonlinear with respect to the unknown nodal degrees of freedom. 
\\
\\ 
As already mentioned in Section \ref{Mass Matrix}, we can compute the third integral in \eqref{PotEnerg_c} related to the dynamic part in different ways:

\begin{enumerate}
	\item First possibility is to evaluate the integral at the centroid $\brm{X}_c$ of the polygon in 2D and of the polyhedra in 3D. The displacement and the acceleration are then evaluated at the centroid and multiplied by the area of the element:
	\begin{equation}
	\int\limits_{\Omega_e} \rho  \ddot{\ten{u}}_\Pi \cdot \ten{u}_\Pi \ d\Omega =
	\rho \ \ddot{\ten{u}}_\Pi (\brm{X}_c) \cdot \ten{u}_\Pi (\brm{X}_c) \ \Omega_e \ .
	\label{eq:vem_dyn_centroid}
	\end{equation}
\item Another possibility is to introduce a sub-triangulation of the polygon or polyhedra and again use one point Gauss integration which yields an evaluation at the centroid $\left. \brm{X}_c\right |_{T}$ of each triangle: 
\begin{equation}
\int\limits_{\Omega_e} \rho  \ddot{\ten{u}}_\Pi \cdot \ten{u}_\Pi \ d\Omega =
\rho \ \sum_T^{nT} \ddot{\ten{u}}_\Pi (\left. \brm{X}_c\right |_{T}) \cdot \ten{u}_\Pi (\left. \brm{X}_c\right |_{T}) \ \Omega_T .
\label{eq:vem_dyn_triangle}
\end{equation}
Since the integral contains quadratic terms of $X$ and $Y$, the two ways introduced above approximate the integral.
\item As a third option, we can compute the integral \eqref{eq:Integral_HTH} {\it exactly}, using the nodal coordinates at the boundary, see \cite{Sing93, Pete13}:
\begin{equation}
\begin{array}{c@{\ }c@{\ }l}
\int\limits_{\Omega_e} 1 \ d \Omega=\frac{1}{2} \sum_{i=1}^{n_V} \left[ y_i x_{i-1} - y_{i-1} x_i \right] \\[2mm] 
\int\limits_{\Omega_e} X \ d \Omega=\frac{1}{6}\sum_{i=1}^{n_V} \left[(y_i x_{i-1} - y_{i-1} x_i)(x_i+x_{i-1})\right] \\[2mm]
\int\limits_{\Omega_e} Y \ d \Omega=\frac{1}{6}\sum_{i=1}^{n_V} \left[(y_i x_{i-1} - y_{i-1} x_i)(y_i+y_{i-1})\right] \\[2mm]
\int\limits_{\Omega_e} X^2 \ d \Omega=\frac{1}{12}\sum_{i=1}^{n_V} \left[(y_i x_{i-1} - y_{i-1} x_i)((x_i + x_{i-1})^2-(x_i x_{i-1}))\right] \\[2mm]
\int\limits_{\Omega_e} Y^2 \ d \Omega=\frac{1}{12}\sum_{i=1}^{n_V} \left[(y_i x_{i-1} - y_{i-1} x_i)((y_i + y_{i-1})^2-(y_i y_{i-1}))\right] \\[2mm]
\int\limits_{\Omega_e} XY \ d \Omega=\frac{1}{24}\sum_{i=1}^{n_V} \left[(y_i x_{i-1} - y_{i-1} x_i)(2 x_{i-1} y_{i-1} + x_{i-1} y_i + x_i y_{i-1} + 2 x_i y_i)\right] ,\\[2mm]
\end{array} 
\label{eq:moments of area}
\end{equation}
where $x_0=x_{n_V}$ and $y_0=y_{n_V}$.
\end{enumerate}

\subsubsection{Stabilization part}
\label{Stabilization part}
The consistency term is computable but yields to a rank deficient {\it stiffness matrix} and thus needs to be stabilized. The idea is to introduce a new positive definite energy $\hat{U}$, with the help of which the stabilization term is redefined, as introduced in \cite{wriggers+reddy+rust+hudobivnik17}:
\begin{equation}
U_{stab} (\boxe{\ten{u}_h} -\UPie) = \hat{U}(\boxe{\ten{u}_h})-\hat{U}(\UPie)
\label{eq:split Ustab with Uhat}
\end{equation}
We further define a stabilization parameter $\beta \in [0,1]$ for the definition of the positive definite energy $\hat{U}$ as:
\begin{equation}
\hat{U}=\beta \ U_c
\label{eq:Uhat}
\end{equation}
Thus the stabilization term takes the form:
\begin{equation}
U_{stab} (\boxe{\ten{u}_h} -\UPie)=\beta U_c(\boxe{\ten{u}_h}) - \beta U_c(\UPie)
\label{eq:split Ustab with Uhat with beta}
\end{equation}
Applying equation \eqref{eq:split Ustab with Uhat with beta} in equation \eqref{eq:Pseudopotential stat+dyn proj+stab}, the final form of the total potential energy function takes the form:

\begin{equation}
U(\ten{u}_h) = \raise3pt
\hbox{$
	\hbox{\scriptsize $n_e$}\atop{\hbox{\LARGEbsf A}\atop
		{\scriptstyle e=1}}$} \Big[ (1-\beta) U_c(\UPie) + \beta U_c(\boxe{\ten{u}_h}) \Big]
\label{eq:Totalenergy}
\end{equation}
The computation of the first term of equation \eqref{eq:Totalenergy} can be done as explained in Section \ref{Consistency part}. The second term $U_c(\boxe{\ten{u}_h})$ needs an approximation. An approach how to compute this part is introduced in \cite{wriggers+reddy+rust+hudobivnik17}. The displacement field is approximated by introducing an internal mesh of 3 noded triangles in 2D or 4 noded tetrahedra in 3D with linear ansatz functions. The nodes of the generated submesh belong to the set of nodes in the virtual element, such that no additional nodes are needed. Based on the triangulated submesh, the displacement gradient is computed. The stabilization term $U_c(\boxe{\ten{u}_h})$ contains both the static $U_c^{stat}(\boxe{\ten{u}_h})$ and dynamic part $U_c^{dyn}(\boxe{\ten{u}_h})$. As the ansatz is linear, the gradient is constant and thus the integral for the static part can be simply evaluated at the centroid $\left. \brm{X}_c\right |_{T}$ of each triangle $n_T$, as sketched in Figure \ref{VEM_EL_Face}. For the dynamic part, the integral can also be evaluated at the centroid of each triangle as shown with elements VEM H2S-I and H2S-II in Section \ref{Examples}. For this case, we would use \eqref{eq:vem_dyn_triangle} but replace the projected quantities with the nodal values for the displacement $\ten{u}_h$ and the acceleration $\ddot {\ten{u}}_h$. The nodal accelerations can be simply computed with the use of \eqref{eq:velocity-Newmark} and \eqref{eq:acceleration-Newmark}.
The stabilization parameter $\beta$ can be chosen freely. For $\beta=1$ the total energy is calculated using only the stabilization part and thus the solution is purely related to the FEM results with three noded triangles in 2D or four noded tetrahedron in 3D. $\beta=0$ yields to a rank deficient stiffness matrix. The choice for the stabilization parameter $\beta$ was analyzed in \cite{hudobivnik2019low,aldakheel+blaz+wriggers18} and it has been shown that the optimal value is in the range $\beta \in [0.2, 0.6]$. For our investigations we choose $\beta=0.4$. Because of the coupling of VEM and FEM, we call this stabilization procedure mixed VEM-FEM-Stabilization.
To obtain the {total} residual vector $\brm{R}_e$ and {total} element tangent matrix $\brm{K}_e$, we compute the first and second derivative of the total energy $U(\ten{u}_h)$ with respect to the global unknowns $\brm{u}_e$ as:

\begin{equation}
\begin{aligned}
\brm{R_e} &= (1-\beta) \underbrace{\left[\brm{R_{c,e}}^{dyn}+\brm{R_{c,e}}^{stat}  \right]}_{\brm{R_{c,e}}}+\beta\underbrace{\left[\brm{R_{stab,e}}^{dyn}+\brm{R_{stab,e}}^{stat}  \right]}_{\brm{R_{stab,e}}}
\\ &=
(1-\beta) \left. \overbrace{\partial_{\brm{u}_e} \left[U_c^{dyn}(\UPie)+U_c^{stat}(\UPie) \right] }^{\partial_{\brm{u}_e} U_c( \UPie)} \right|_{\ddot {\brm{u}}_e = const.} \\
&+ \beta \left. \overbrace{  \partial_{\brm{u}_e} \left[ U_c^{dyn}(\boxe{\ten{u}_h})+U_c^{stat}(\boxe{\ten{u}_h}) \right]} ^{\partial_{\brm{u}_e} U_c(\boxe{\ten{u}_h})} \right|_{\ddot {\brm{u}}_e = const.} 
\end{aligned}
\label{eq:Residualvector}
\end{equation}
\\
\begin{equation}
\begin{aligned}
\brm{K_e} &= (1-\beta) \underbrace{\left[\brm{K_{c,e}}^{dyn}+\brm{K_{c,e}}^{stat}  \right]}_{\brm{K_{c,e}}}+\beta\underbrace{\left[\brm{K_{stab,e}}^{dyn}+\brm{K_{stab,e}}^{stat}  \right]}_{\brm{K_{stab,e}}}
\\ &=
(1-\beta) \dfrac{\partial \brm{R_{c,e}}}{ \partial \brm{u}_e} + \beta  \dfrac{\partial \brm{R_{stab,e}}}{ \partial \brm{u}_e} 
\end{aligned}
\label{eq:Tangentmatrix}
\end{equation}

For the separate computation of the mass- and the stiffness matrix, we split $\beta$ into $\beta^{stat}$ and $\beta^{dyn}$. Equation \eqref{eq:Residualvector} and \eqref{eq:Tangentmatrix} become:

\begin{equation}
\begin{aligned}
\brm{R_e} &= (1-\beta^{stat}) \dfrac{\partial U_c^{stat}(\UPie) }{ \partial \brm{u}_e} + \left. (1-\beta^{dyn})   \dfrac{\partial U_c^{dyn}(\UPie)}{ \partial \brm{u}_e} \right|_{\ddot {\brm{u}}_e = const.} \\
&+ \beta^{stat} \dfrac{\partial U_c^{stat}(\boxe{\ten{u}_h}) }{ \partial \brm{u}_e} + \left. \beta^{dyn}   \dfrac{\partial U_c^{dyn}(\UPie)}{ \partial \brm{u}_e} \right|_{\ddot {\brm{u}}_e = const.}
\end{aligned}
\label{eq:New Residualvector}
\end{equation}

\begin{equation}
\begin{aligned}
\brm{K_e} &= (1-\beta^{stat}) \brm{K_{c,e}}^{stat} + (1-\beta^{dyn}) \brm{K_{c,e}}^{dyn}+ \beta^{stat} \brm{K_{stab,e}}^{stat} + \beta^{dyn} \brm{K_{stab,e}}^{dyn} 
\end{aligned}
\label{eq:New Tangentmatrix}
\end{equation}
Here, $\beta^{stat}$ and $\beta^{dyn}$ are the stabilization parameters for the static and dynamic part, where $\beta^{stat}, \beta^{dyn} \in [0,1]$.
For the computation of the derivatives in equation \eqref{eq:New Residualvector} and \eqref{eq:New Tangentmatrix}, the Mathematica Package \AceGenNamec \ is used, see \cite{KoWr16}.
\section{Numerical Examples}
\label{Examples}
In this section we demonstrate the performance of the presented 3D virtual element formulation for dynamic problems at finite deformations. For comparison purposes results of the standard finite element method (FEM) are also included. The material parameters used in this work are the same for all examples and are given in Table \ref{table:MP1}.

\begin{table}[!htb]
	\begin{center}
		\caption{Material parameters used for the numerical examples}
		\label{table:MP1}
		\begin{tabular}{l l l r c}\hline
			No. &Parameter & Label & Value & Unit \\\hline
			1 & Elastic modulus & $E$ & $210$ & $kN/mm^2$\\
			2 & Poisson ratio & $ \nu $ & $0.3$ & --\\
			3 & Density & $\rho$  & $0.0027$ & $g/mm^3$\\\hline
		\end{tabular}
	\end{center}
\end{table}

In this contribution, the following mesh types with first order virtual element discretizations are introduced. Herein, the different ways to evaluate the integral in \eqref{eq:Integral_HTH} are employed in the following element types:

\begin{itemize}
	\item VEM Q2S: A regular 2D mesh with 8 noded quadrilateral elements. The argument of the integral is evaluated at the {\it centroid} of the polygon and is multiplied by the area of the element. $\beta^{dyn}=0$, which means that the mass-matrix is computed using only the projection part. This represents the classical way as introduced in \eqref{eq:vem_dyn_centroid}.
	\item VEM Q2S BI: A regular 2D mesh with 8 noded quadrilateral elements. The argument of the integral is computed {\it exactly} on the boundary with the moments of area \eqref{eq:moments of area}. $\beta^{dyn}=0$, which means that the mass-matrix is computed using only the projection part.
	\item VEM VO: 2D/3D Voronoi cell mesh with arbitrary number of element nodes. The argument of the integral is evaluated at the centroid of the polygon and is multiplied by the area of the element. $\beta^{dyn}=0$, which means that the mass-matrix is computed using only the projection part.
	\item VEM VO BI: 2D/3D Voronoi cell mesh with arbitrary number of element nodes. The argument of the integral is computed {\it exactly} on the boundary with the moments of area \eqref{eq:moments of area}. $\beta^{dyn}=0$, which means that the mass-matrix is computed using only the projection part.
	\item VEM H2S: A regular 3D mesh with 20 noded hexahedral elements. The argument of the integral is evaluated at the centroid of the polygon and is multiplied by the area of the element. $\beta^{dyn}=0$, which means that the mass-matrix is computed using only the projection part.
	\item VEM H2S-I: A regular 3D mesh with 20 noded hexahedral elements. The polyhedra is {\it subdivided} into tetrahedrons and the argument of the integral is evaluated at the {\it centroid of each tetrahedron} and is multiplied by the volume of the element. $\beta^{dyn}=0$, which means that the mass-matrix is computed using only the projection part. This procedure is same as in \eqref{eq:vem_dyn_triangle}.\\
	\item VEM H2S-II: A regular 3D mesh with 20 noded hexahedral elements. The polyhedra is subdivided into tetrahedrons and the argument of the integral is evaluated at the centroid of each tetrahedron and is multiplied by the volume of the element. $\beta^{dyn}=1$, which means that the mass-matrix is computed using only the stabilization part.\\
\end{itemize}

For a representative comparison, the following finite element formulations were selected:

\begin{itemize}
	\item FEM T1: A regular 2D mesh with 3 noded triangular first order finite elements.
	\item FEM Q1: A regular 2D mesh with 4 noded quadliteral first order finite elements.
	\item FEM Q2: A regular 2D mesh with {8} noded quadliteral second order finite elements.
	\item FEM H1: A regular 3D mesh with 8 noded first order finite elements.
	\item FEM H2: A regular 3D mesh with 27 noded second order finite elements.
\end{itemize}

For the stabilization parameter of the {\it static part} we choose $\beta^{stat}=0.4$ for all the simulations. Unless otherwise specified, all the results for the dynamic part are obtained with $\beta^{dyn}=0$. 
In this case the mass-matrix is computed using only the projection part having no stabilization. Therefore, equations \eqref{eq:New Residualvector} and \eqref{eq:New Tangentmatrix} simplify to:

\begin{equation}
\begin{aligned}
\brm{R_e} &= (1-\beta^{stat}) \dfrac{\partial U_c^{stat}(\UPie) }{ \partial \brm{u}_e} + \left.   \dfrac{\partial U_c^{dyn}(\UPie)}{ \partial \brm{u}_e} \right|_{\ddot {\brm{u}}_e = const.} + \beta^{stat} \dfrac{\partial U_c^{stat}(\boxe{\ten{u}_h}) }{ \partial \brm{u}_e}
\end{aligned}
\label{eq:New2 Residualvector}
\end{equation}

\begin{equation}
\begin{aligned}
\brm{K_e} &= (1-\beta^{stat}) \brm{K_{c,e}}^{stat} + \brm{K_{c,e}}^{dyn}+ \beta^{stat} \brm{K_{stab,e}}^{stat}
\end{aligned}
\label{eq:New2 Tangentmatrix}
\end{equation}

\begin{figure}[b]
	\centering
	\scalebox{0.9}{\includegraphics{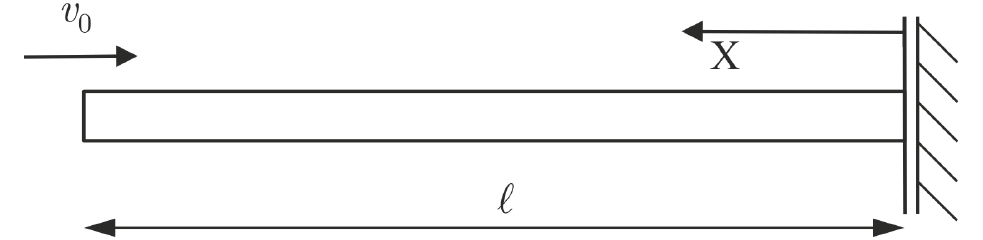}} 
	\caption{2D Example - Wave propagation in longitudinal beams (Boundary value problem).}
	\label{2D Example Longitudinal - BVP}
\end{figure}
\begin{figure}[!t]
	\begin{minipage}{0.5\textwidth}
		\hspace{30mm}
		\centering
		\begin{subfigure}[c]{0.5\textwidth}
			\scalebox{1}{\includegraphics{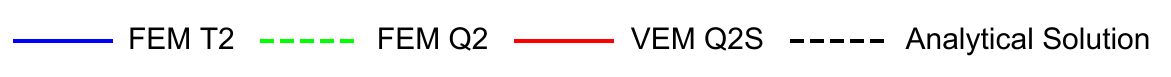}}
			\label{2D Example Longitudinal - Legend}
		\end{subfigure}
	\end{minipage}
	
	\begin{subfigure}[c]{0.5\textwidth}
		\scalebox{0.62}{\includegraphics{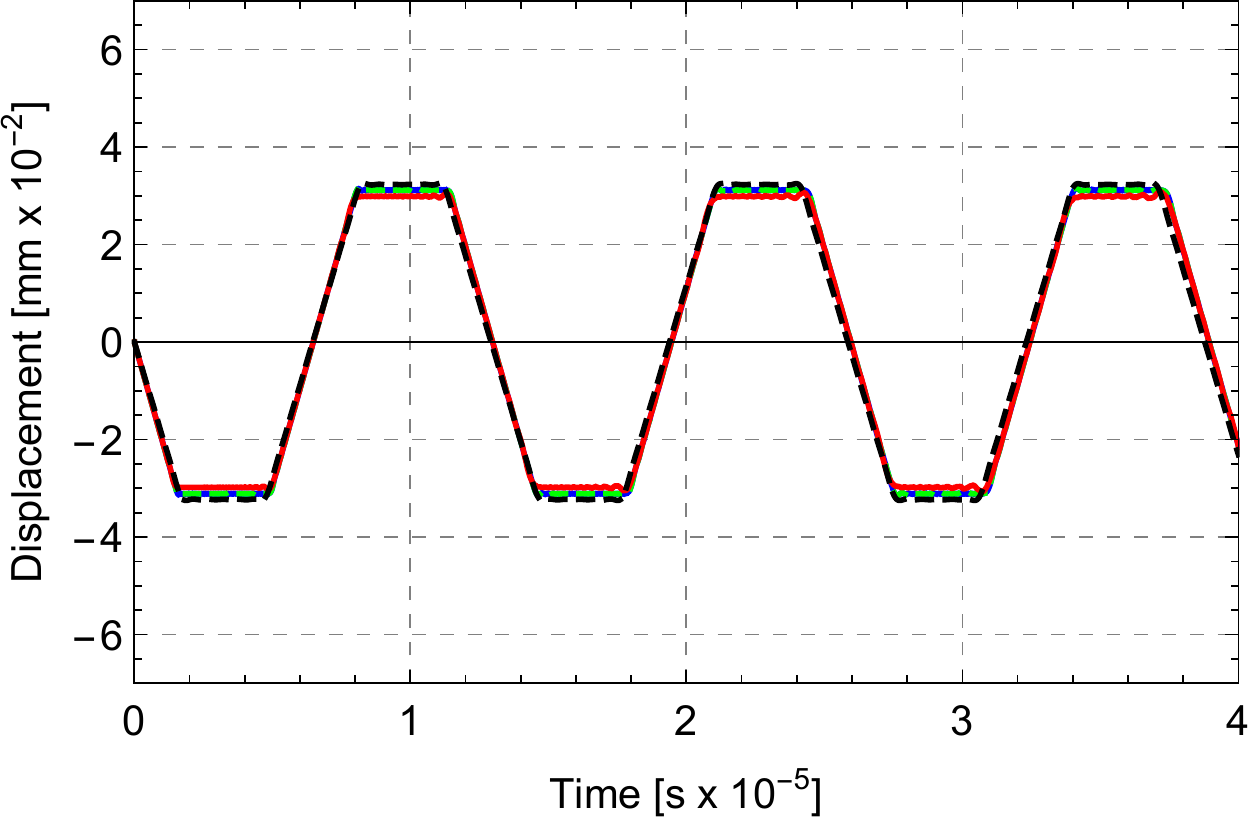}} 
		\subcaption{Response at $x=\ell/2$}
		\label{2D Example Longitudinal - Comparison 1}
	\end{subfigure}
	\begin{subfigure}[c]{0.5\textwidth}
		\scalebox{0.62}{\includegraphics{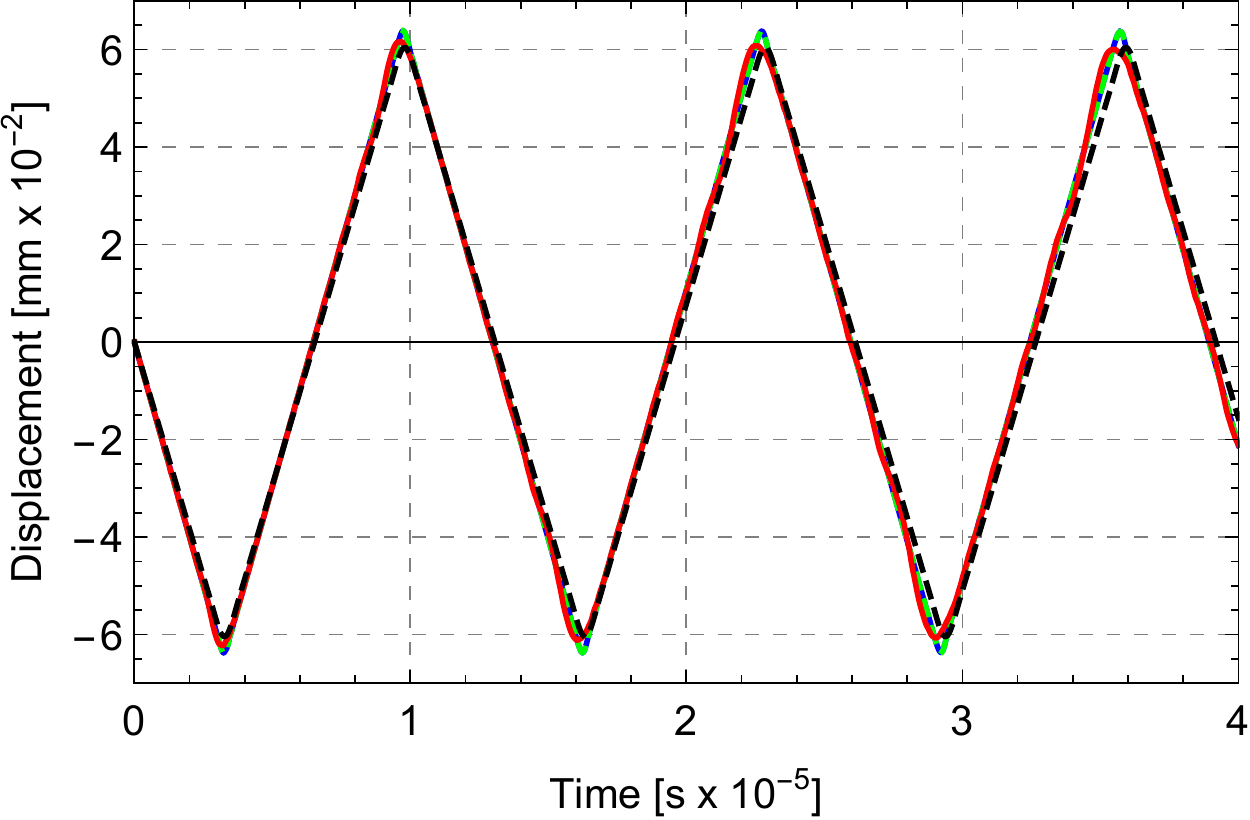}} 
		\subcaption{Response at $x=\ell$}
		\label{2D Example Longitudinal - Comparison 2}
	\end{subfigure}
	\caption{ Displacement over time response for 2D Example - Wave propagation in longitudinal beams.}
\end{figure}

\subsection{2D Boundary value problems}
\subsubsection{Wave propagation in longitudinal beams}

The first model problem is concerned with analyzing the wave propagation in longitudinal beams. The geometric setup and the loading conditions of the specimen are depicted in Figure \ref{2D Example Longitudinal - BVP}. The height of the beam is chosen to be $h=0.3 \ mm$ and the length $\ell=30 \ mm$, where the degrees of freedom are fixed in longitudinal direction on the right side. The height of the specimen is much smaller than the length thus this problem is considered to be a one dimensional problem. The initial velocity is set to $v_0=20 \ m/s$. After a certain time, we observe the wave propagation through the elastic body. Next, the virtual element method will be compared with the finite element method and the analytical solution. The analytical solution in \eqref{eq:analyticalsolution1} is obtained by solving the wave equation \eqref{eq:Waveequation}. Figure \ref{2D Example Longitudinal - Comparison 1} and  \ref{2D Example Longitudinal - Comparison 2} illustrate the displacement field over time for different VEM and FEM formulations and compared with analytical results as well. 
The FEM results are computed for $4 \times 200$  elements, where the virtual element results are computed for $4 \times 100$ elements. We observe a good agreement of VEM compared with FEM solution and the analytical solution. In terms of the period and the amplitude of the wave, the virtual elements shows results that are close to the analytical solution
\begin{equation}
\dfrac{\partial^2 u}{\partial t^2}=c^2 \dfrac{\partial^2 u}{\partial x^2}
\quad \textrm{where} \quad 
c=\sqrt{\frac{E}{\rho}}
\label{eq:Waveequation}
\end{equation}
\begin{equation}
u(x,t) = \sum_{n=0}^{\infty} \frac{2 v_0 c}{\ell {\omega_n}^2} sin \left(\frac{w_n x}{c} \right) sin(w_n t) 
\quad \textrm{with} \quad
w_n = \frac{1}{0,95} \frac{(2n+1) \pi c}{2 \ell} \ .
\label{eq:analyticalsolution1}
\end{equation}
In VEM Q2S the integral for the dynamic part in equation \eqref{eq:Integral_HTH} is evaluated at the centroid of the element, hence this method seems to be sufficient.

\subsubsection{Transversal Beam Vibration}

\begin{figure}[!t]
	\begin{subfigure}[c]{0.5\textwidth}
		\scalebox{1.0}{\includegraphics{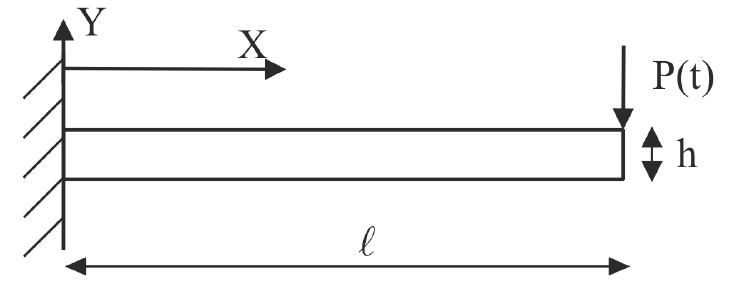}} 
		\subcaption{}
		\label{2D Example Transversal - BVP}
	\end{subfigure}
	\begin{subfigure}[c]{0.5\textwidth}
		\scalebox{0.65}{\includegraphics{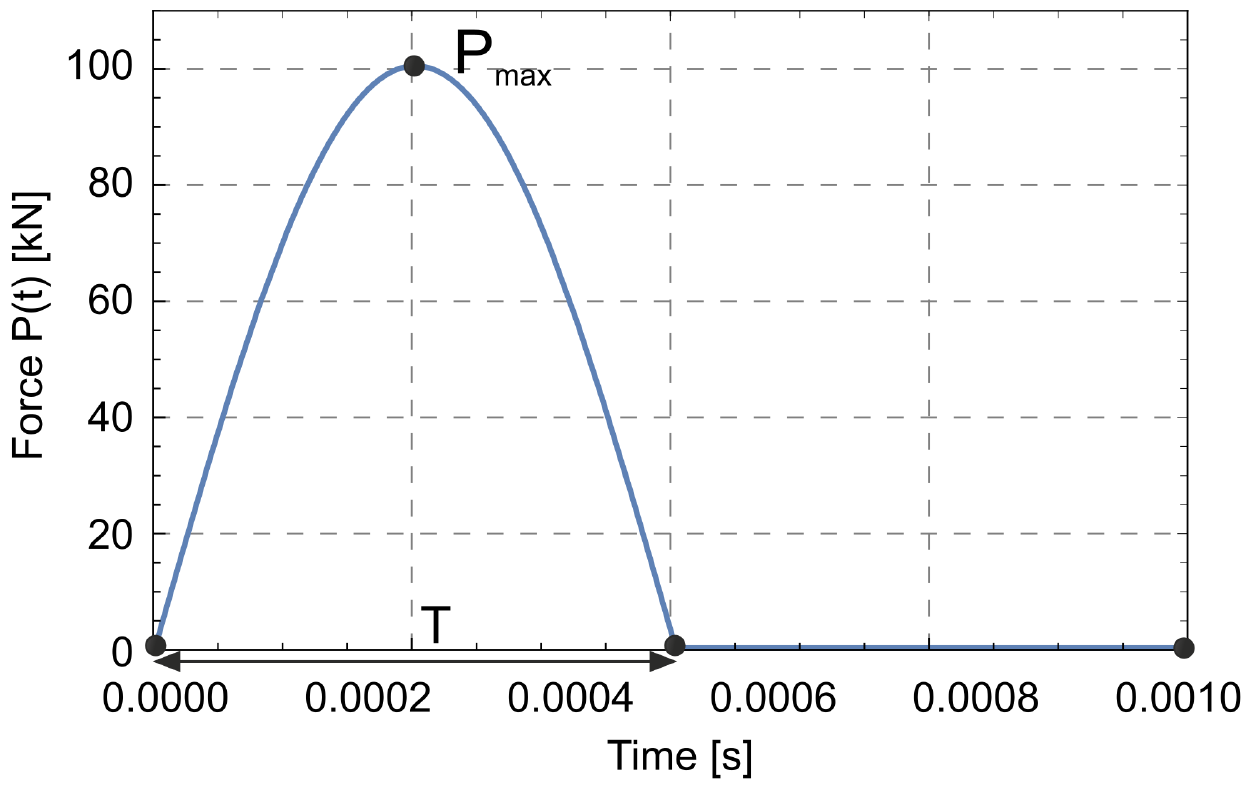}} 
		\subcaption{}
		\label{2D Example Transversal - Force}
	\end{subfigure}
	\caption{2D Example - Transversal beam vibration. Boundary value problem in \textbf{(a)} and applied force in \textbf{(b)}.}
\end{figure}
\begin{figure}[!t]
	\begin{subfigure}[t]{1\textwidth}
		\centering
		\scalebox{1.0}{\includegraphics{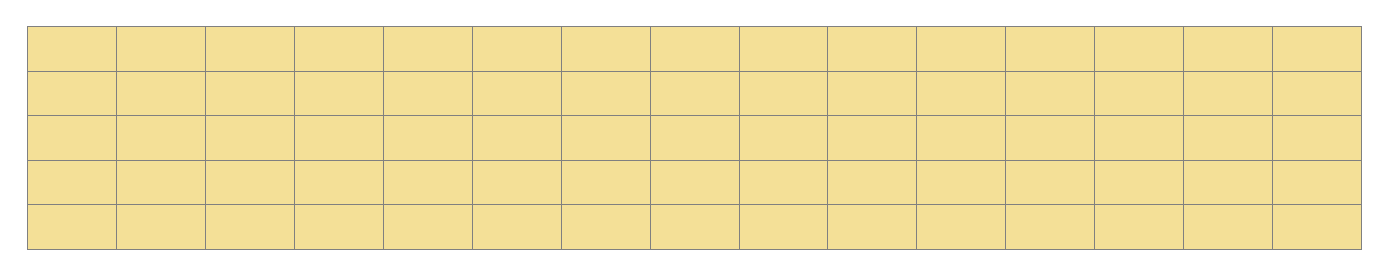}}
		\subcaption{}
		\label{2D Example Transversal Beam Vibration - VEM Q2S Mesh}
	\end{subfigure}
	\begin{subfigure}[t]{1\textwidth}
		\centering
		\scalebox{1.1}{\includegraphics{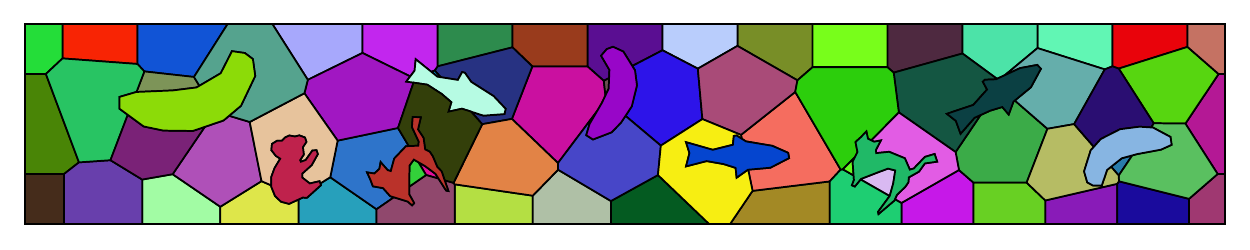}}		 
		\subcaption{}
		\label{2D Example Transversal Beam Vibration - Animal Mesh}		
	\end{subfigure}
	\begin{subfigure}[t]{1\textwidth}
	\centering
		\scalebox{1.1}{\includegraphics{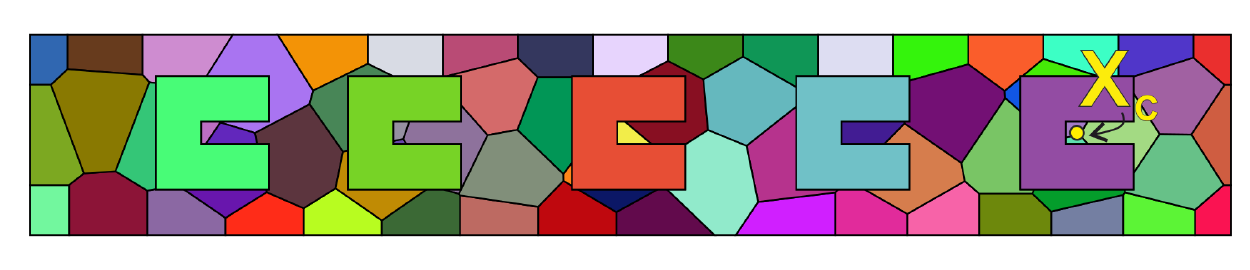}} 
		\subcaption{}
		\label{2D Example Transversal Beam Vibration - C Mesh}
	\end{subfigure}
	\caption{ 2D Example - Transversal Beam Vibration. VEM Q2S-Mesh (\textbf{a}), VEM Animal-Mesh (\textbf{b}) and C-Mesh  (\textbf{c}).}
	\label{2D Example Transversal Beam Vibration - Mesh Types}
\end{figure}

\begin{figure}[!t]
	\begin{minipage}{0.5\textwidth}
		\hspace{35mm}
		\centering
		\begin{subfigure}[c]{0.2\textwidth}
			\scalebox{1}{\includegraphics{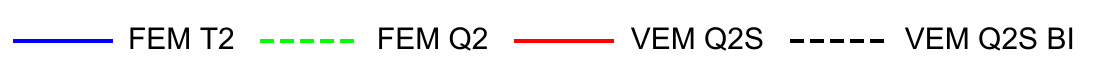}}
			\label{2D Example Transversal - Legend}
		\end{subfigure}
	\end{minipage}
	
	\begin{subfigure}[c]{0.5\textwidth}
		\scalebox{0.62}{\includegraphics{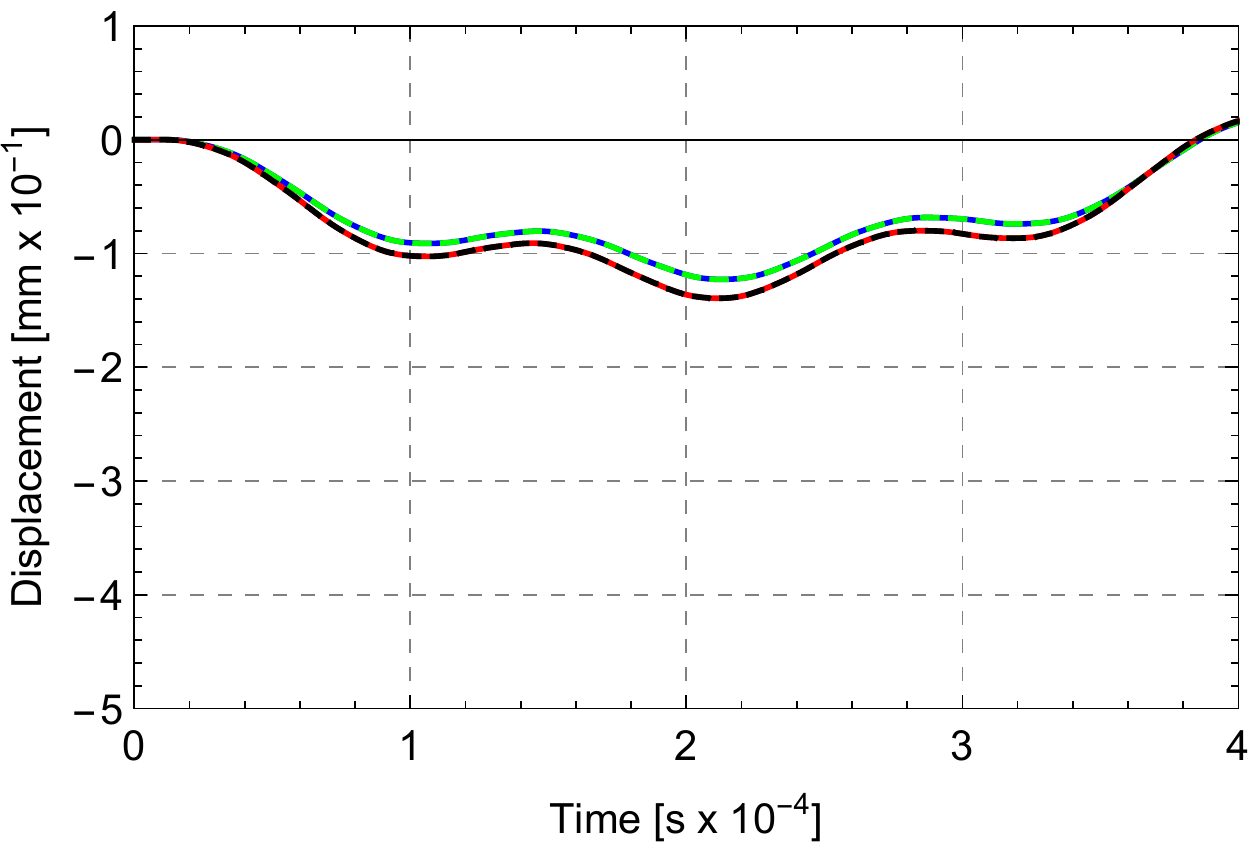}} 
		\subcaption{Response at $x=\ell/2$}
		\label{2D Example Transversal - Comparison 1}
	\end{subfigure}
	\begin{subfigure}[c]{0.5\textwidth}
		\scalebox{0.62}{\includegraphics{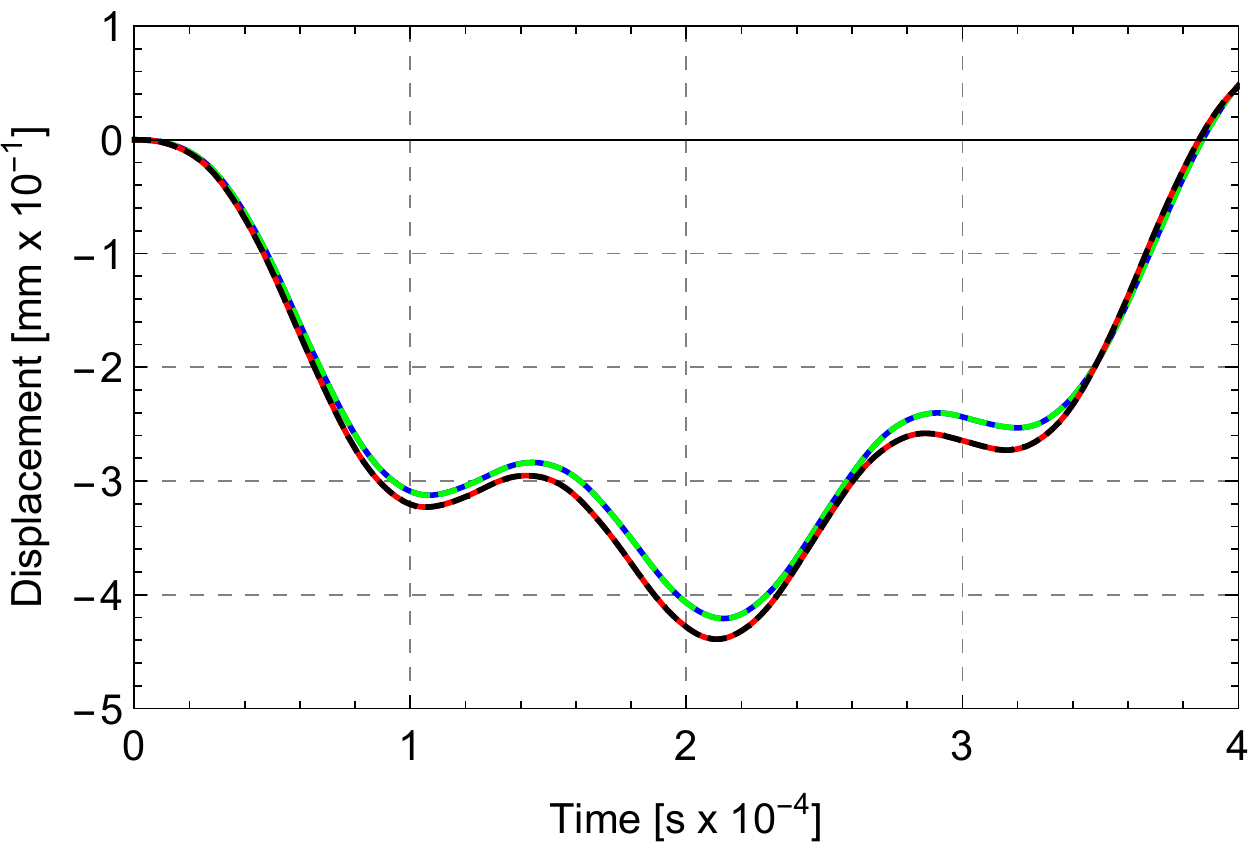}} 
		\subcaption{Response at $x=\ell$}
		\label{2D Example Transversal - Comparison 2}
	\end{subfigure}
	\caption{ Displacement over time response for 2D Example - Transversal Beam Vibration}
\end{figure}
\begin{figure}[!t]
	\begin{minipage}{0.5\textwidth}
		\hspace{35mm}
		\centering
		\begin{subfigure}[c]{0.2\textwidth}
			\scalebox{1}{\includegraphics{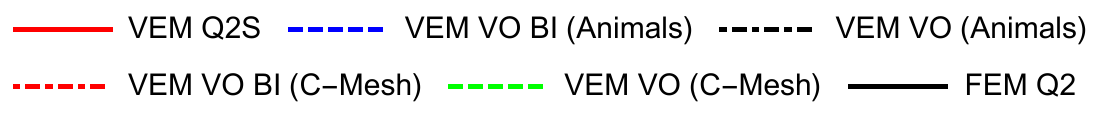}}
			\label{2D Example Transversal - Animal Mesh Legend}
		\end{subfigure}
	\end{minipage}
	
	\begin{subfigure}[c]{0.5\textwidth}
		\scalebox{0.62}{\includegraphics{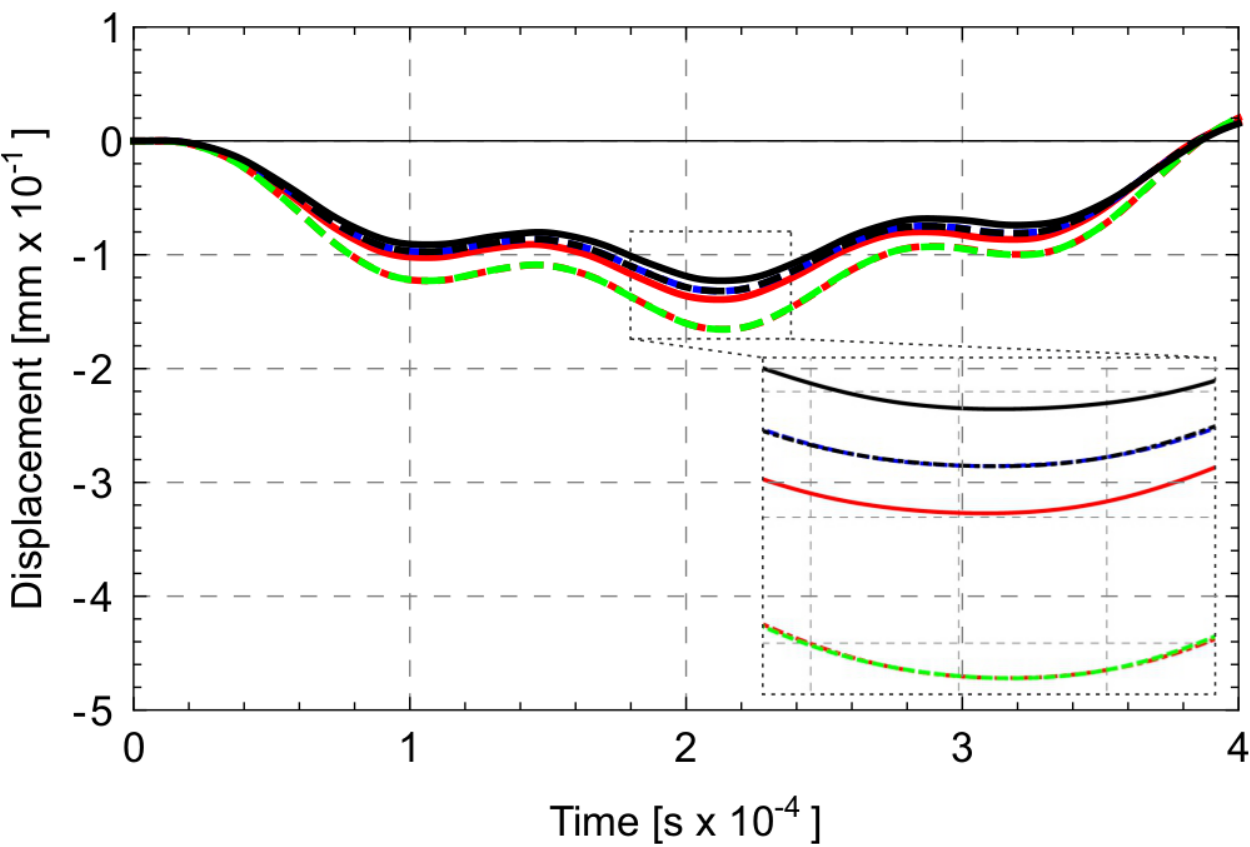}} 
		\subcaption{Response at $x=\ell/2$}
		\label{2D Example Transversal - Animal Comparison 1}
	\end{subfigure}
	\begin{subfigure}[c]{0.5\textwidth}
		\scalebox{0.62}{\includegraphics{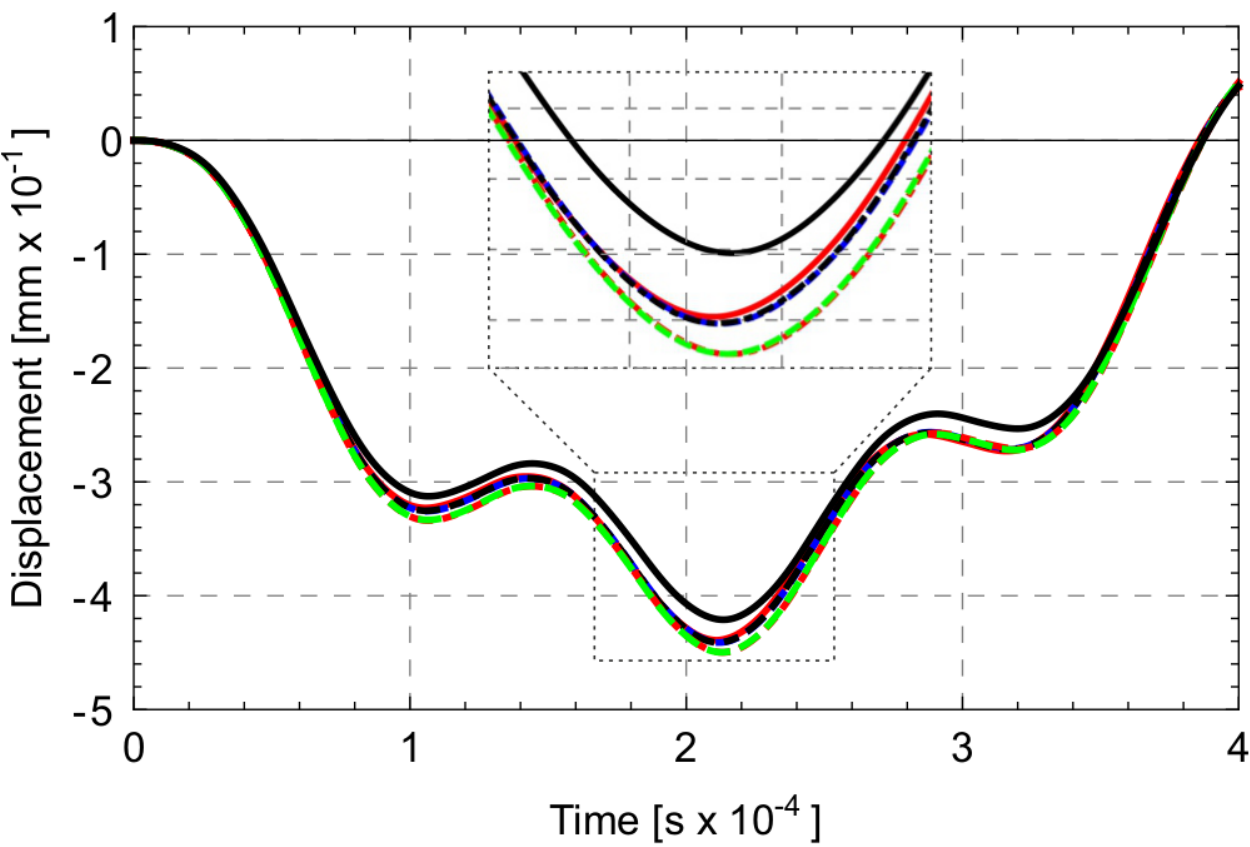}} 
		\subcaption{Response at $x=\ell$}
		\label{2D Example Transversal - Animal Comparison 2}
	\end{subfigure}
	\caption{ Displacement over time response for 2D Example - Transversal Beam Vibration}
\end{figure}
The second benchmark test is concerned with analyzing the transversal vibration in beams. The geometric setup and the loading conditions of the beam are depicted in Figure \ref{2D Example Transversal - BVP}. The length of the bar is set to $\ell=30 \ mm$ and the hight is $h=5 \ mm$. The force is applied transversal at the end of the specimen as shown in Figure \ref{2D Example Transversal - Force}. The temporal course of the force is given by a {\it half sine}, where the maximum of the force is set to $P_{max} = 100 \ kN$. The period $T$ of the applied force is adjusted to the bending stiffness of the beam and defined as

\begin{equation}
T=\frac{3.5156}{2 \pi \ell^2} {\sqrt{\frac{12 \rho}{ E bh^3}}}
\label{eq:Timeperiod Beam}
\end{equation}

In order to analyze the position effect of the element centroid on evaluating the integral of the dynamic part, we used different type of meshes which can be seen in Figure \ref{2D Example Transversal Beam Vibration - Mesh Types}. The "animal" mesh (Figure \ref{2D Example Transversal Beam Vibration - Animal Mesh}) includes non convex elements. To see the effect of using non-convex elements where the centroid of the element is outside of the element domain, we use a special mesh with elements like C's, where the centroid of the element is outside of the element domain (Figure \ref{2D Example Transversal Beam Vibration - C Mesh}). 

Figure \ref{2D Example Transversal - Comparison 1} and \ref{2D Example Transversal - Comparison 2} show the displacement over time response in the center at $x=\ell/2$ and at the end of the beam at $x=\ell$. The finite element solution is computed for $1000$ elements, whereas VEM results are obtained with $100$ virtual elements. The comparison of the virtual elements Q2S and Q2S BI shows that it makes no difference for regular shaped meshes if the integral in equation \eqref{eq:Integral_HTH} is evaluated approximately on the centroid of the element or exactly on the boundary using the moments of area. Furthermore, we can see that the displacements in the center of the beam are slightly higher than the finite element results. However the period fits very well compared with FEM results. In general the virtual element results are in a good agreement with the compared finite element results. 

The comparison of the different meshes shows, that the C-mesh yields a higher deflection, compared to the other results. Nevertheless qualitatively the shape of the displacement over time response fits very well the finite element Q2 results and the virtual element Q2S results. Again, the evaluation of the integral at the centroid of the element compared to computing the integral at the boundary exactly using the moments of area does not affect the results.
\begin{figure}[!b]
	\begin{subfigure}[t]{0.28\textwidth}
		\scalebox{0.35}{\includegraphics{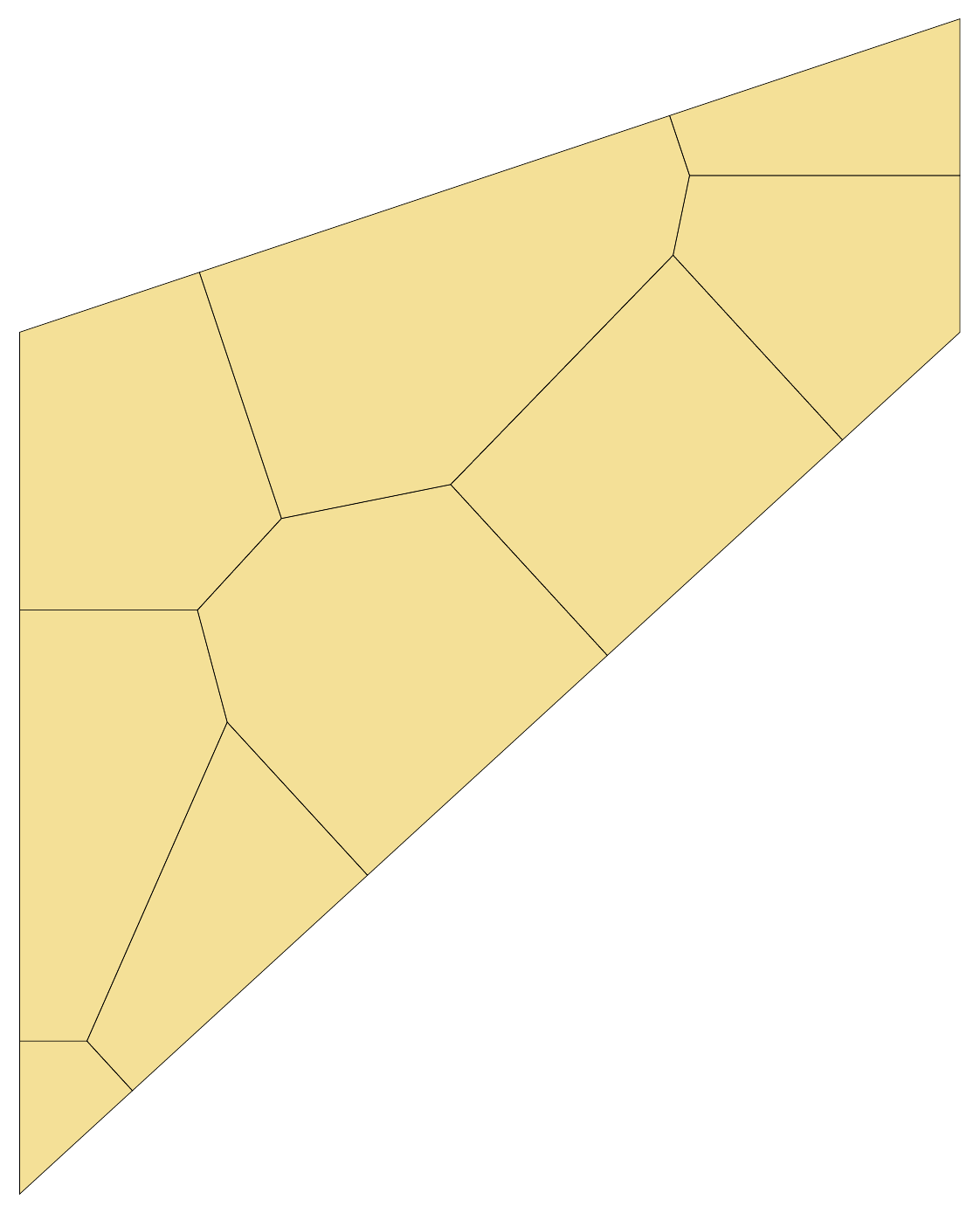}}
		\subcaption{}
		\label{2D Example Cooksmembrane - VEMVoronoi Mesh}
	\end{subfigure}
	\begin{subfigure}[t]{0.45\textwidth}
		\scalebox{1.0}{\includegraphics{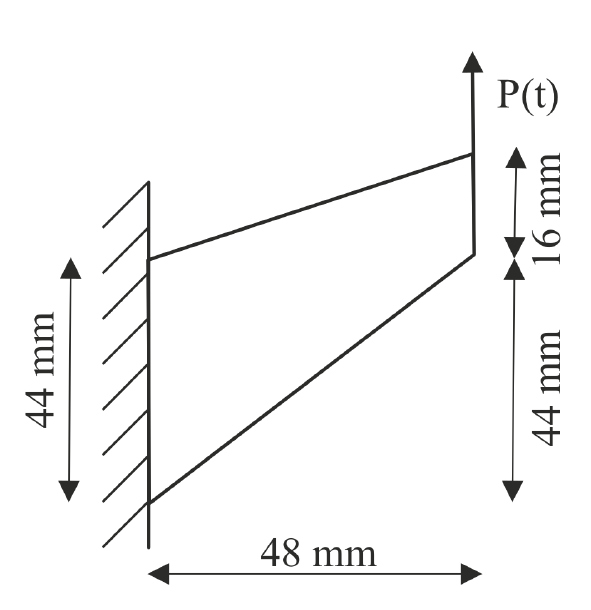}} 
		\subcaption{}
		\label{2D Example Cooksmembrane - BVP}
	\end{subfigure}
	\begin{subfigure}[t]{0.25\textwidth}
		\scalebox{0.35}{\includegraphics{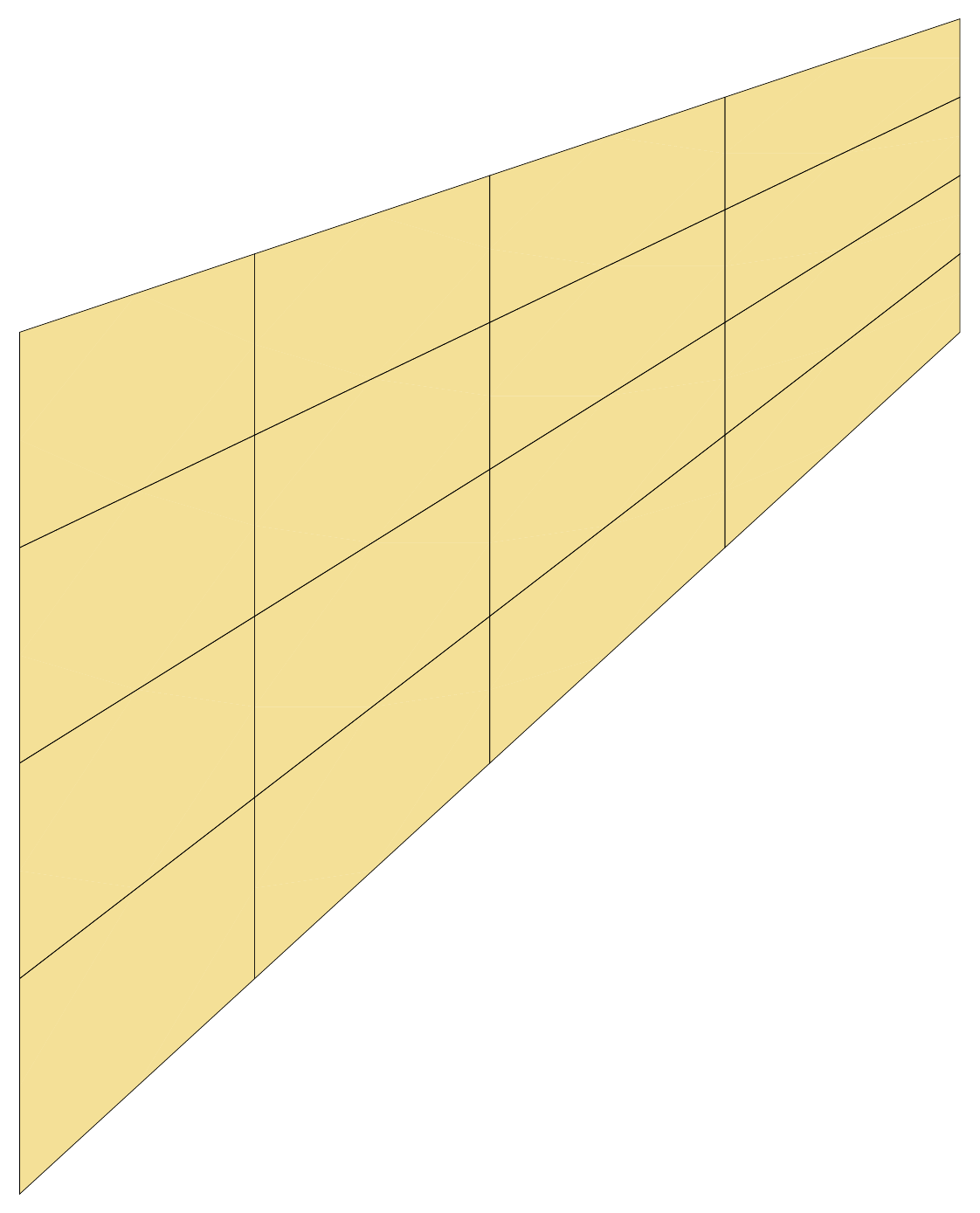}} 
		\subcaption{}
		\label{2D Example Cooksmembrane - VEMQ2S Mesh}
	\end{subfigure}\\[3mm]
	\caption{2D Example - Cook's membrane. VEM Voronoi Mesh (\textbf{a}), Boundary Value Problem (\textbf{b}) and VEM Q2S Mesh (\textbf{c}).}
\end{figure}
\begin{figure}[!t]
	\begin{subfigure}[t]{0.32\textwidth}
		\scalebox{0.6}{\includegraphics{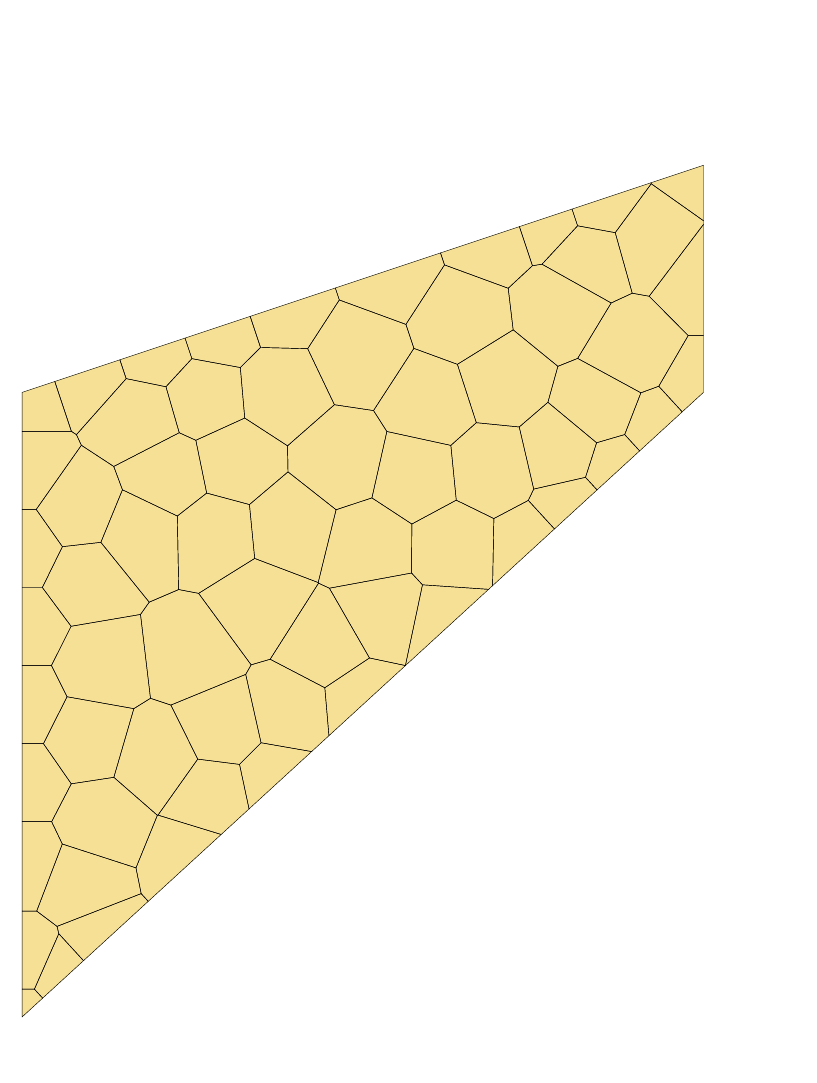}}
		\subcaption{t=0s}
	\end{subfigure}
	\begin{subfigure}[t]{0.32\textwidth}
		\scalebox{0.6}{\includegraphics{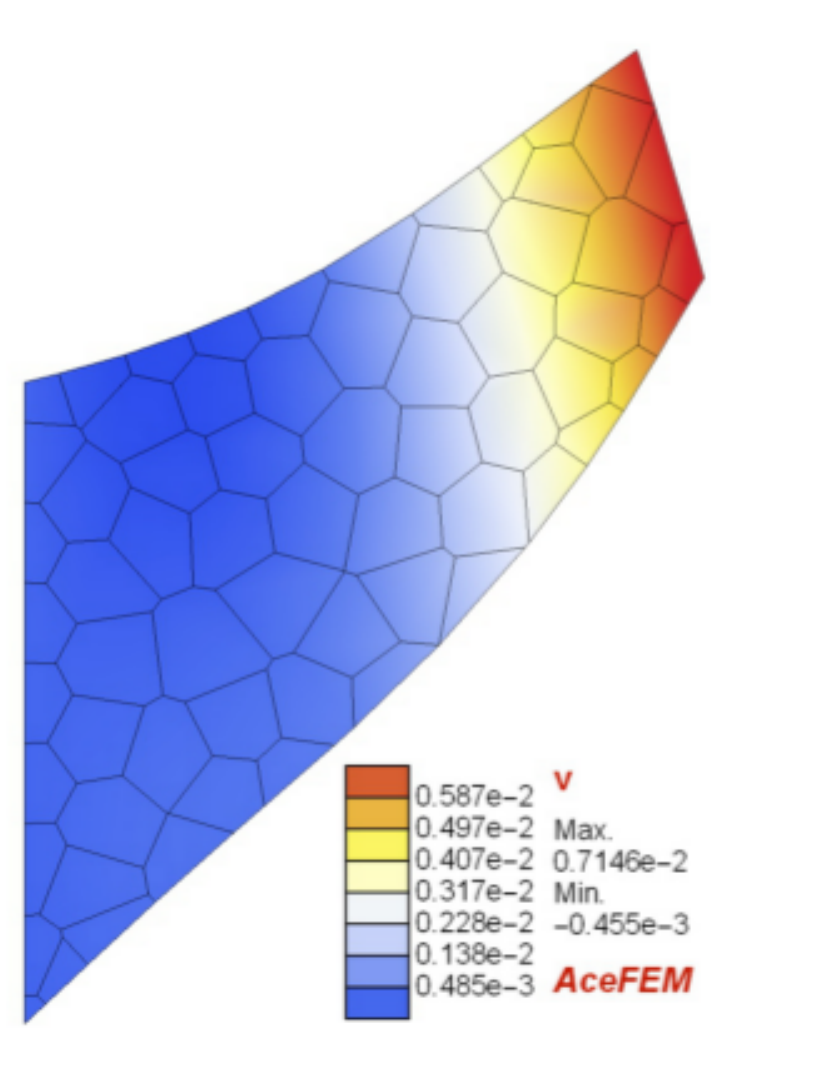}} 
		\subcaption{t=0.0001s}
	\end{subfigure}
	\begin{subfigure}[t]{0.32\textwidth}
		\scalebox{0.6}{\includegraphics{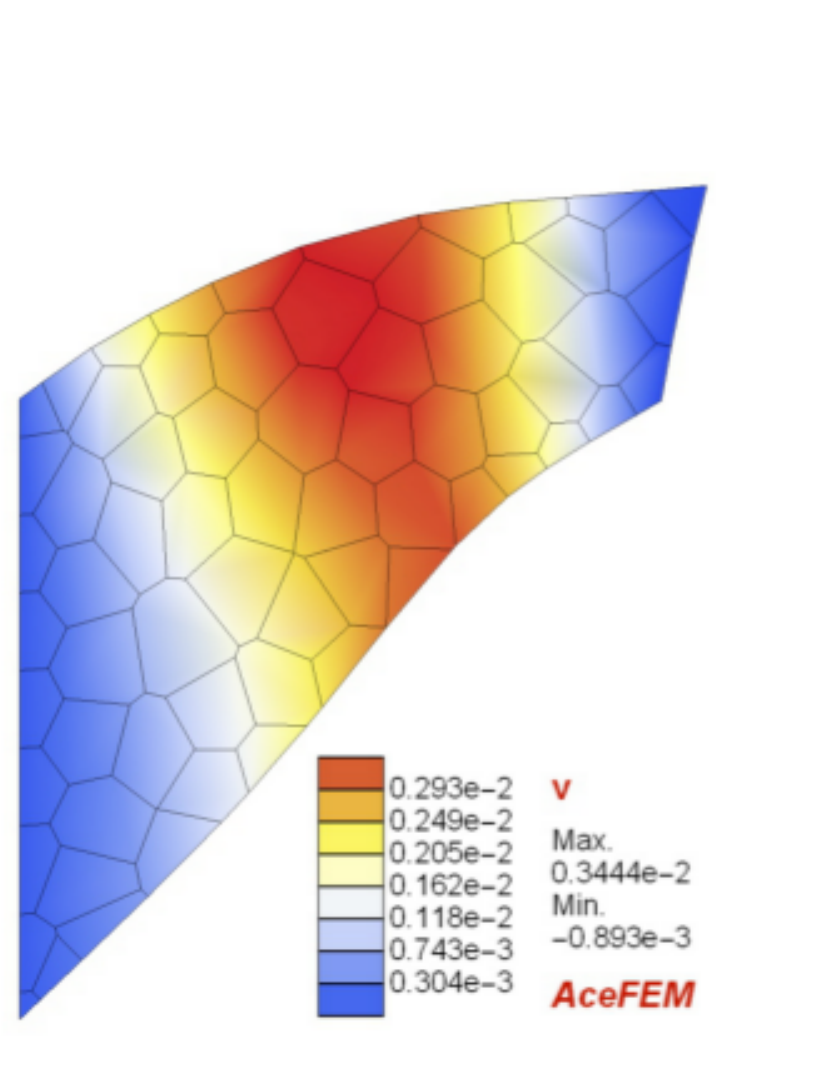}} 
		\subcaption{t=0.0002s}
	\end{subfigure}
	\begin{subfigure}[t]{0.32\textwidth}
		\scalebox{0.6}{\includegraphics{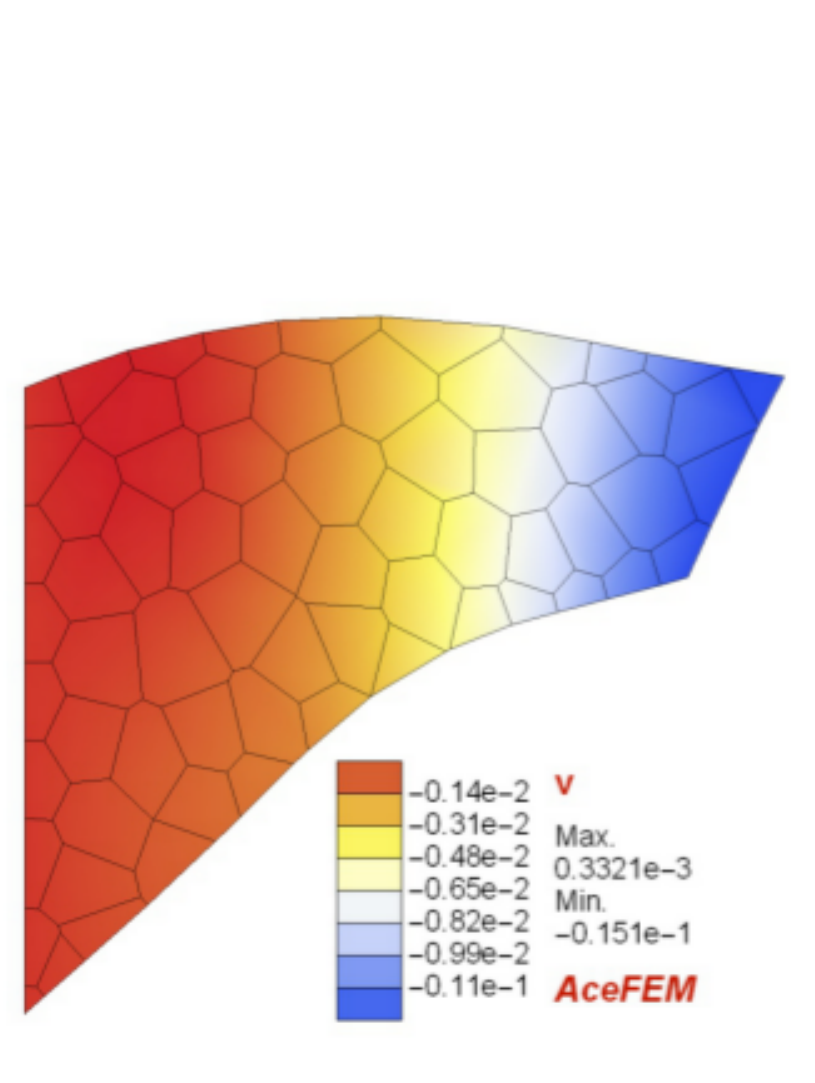}} 
		\subcaption{t=0.00035s}
	\end{subfigure}
	\begin{subfigure}[t]{0.32\textwidth}
		\scalebox{0.6}{\includegraphics{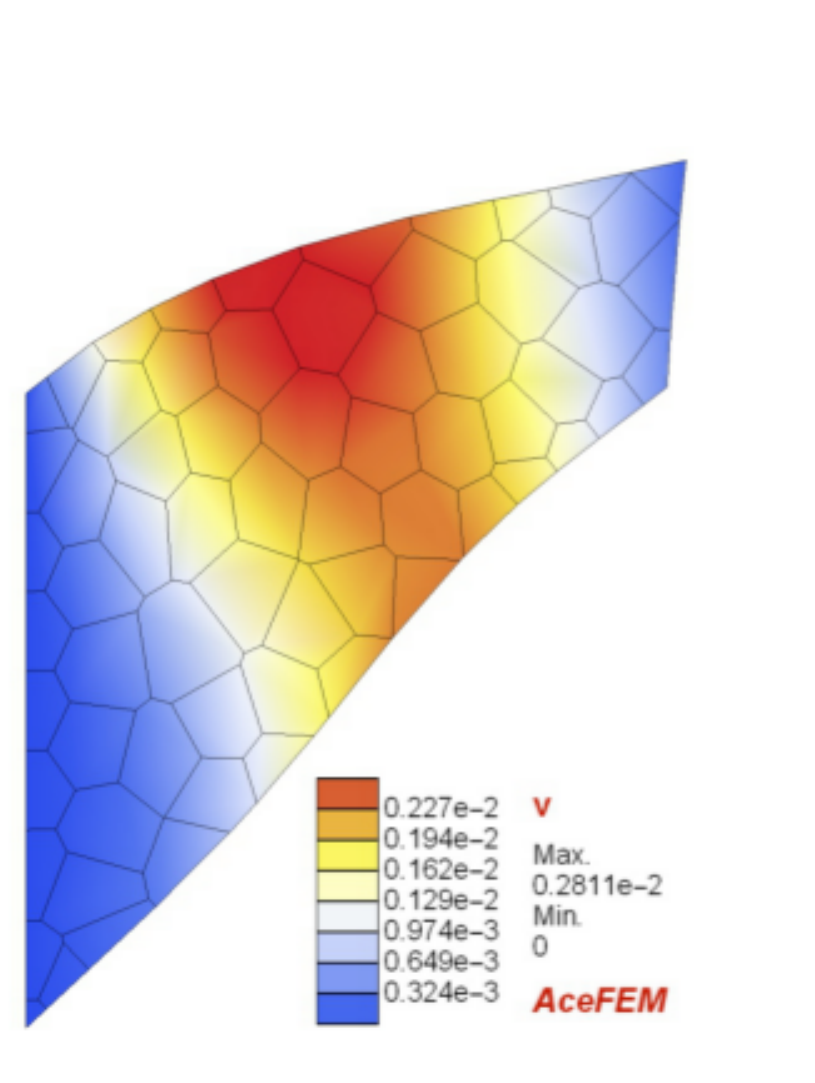}} 
		\subcaption{t=0.00055s}
	\end{subfigure}
	\begin{subfigure}[t]{0.32\textwidth}
		\scalebox{0.6}{\includegraphics{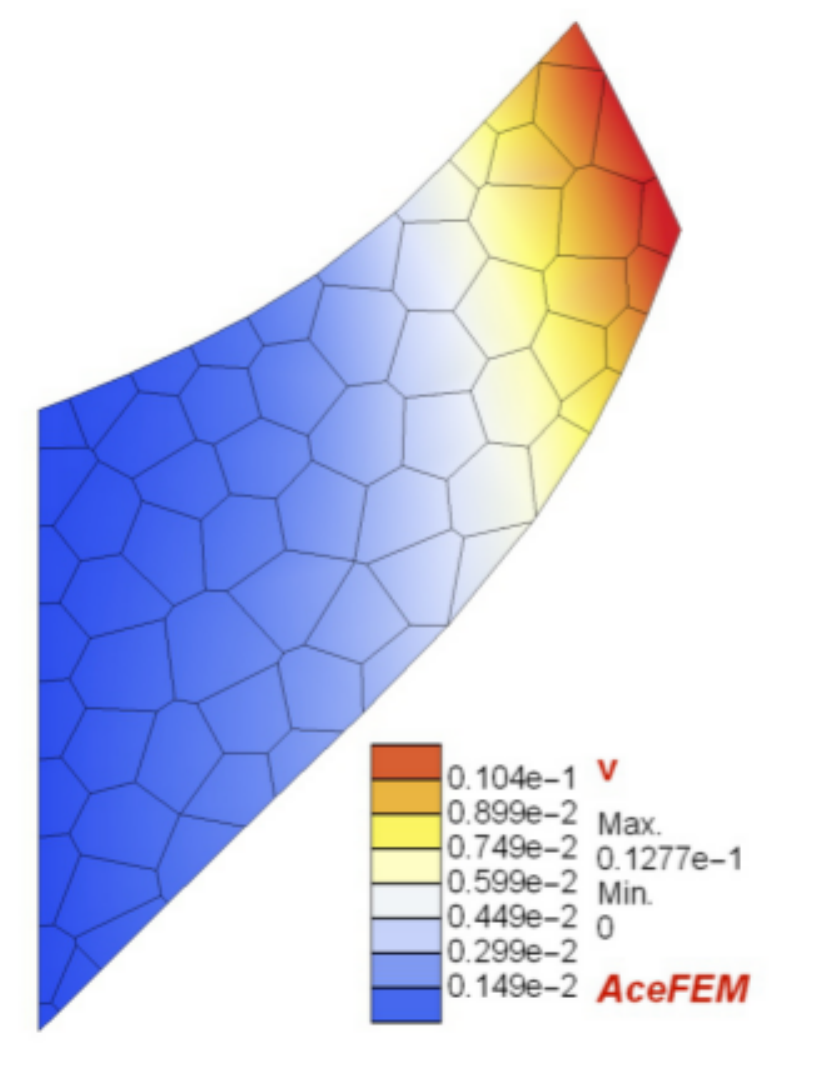}} 
		\subcaption{t=0.00065s}
	\end{subfigure}
	\\[4mm]
\caption{2D Example - Cook's membrane. {Evolution of the vertical displacement for different deformation states.}}
\label{CM-DM}
\end{figure}

\subsubsection{Cook's membrane problem}
The next example is the Cook's membrane problem in 2D. Here as well the virtual element performance will be compared with the finite element results. The geometrical setup and boundary conditions are demonstrated in Figure \ref{2D Example Cooksmembrane - BVP}. In this test a force driven scenario is applied at the right edge as a line load as depicted in Figure \ref{2D Example Cooksmembrane - BVP}. The force is applied as shown in \ref{2D Example Transversal - Force} with $P_{max} = 10 000  \ kN/mm$. VEM VO mesh and regular VEM Q2S mesh are also plotted in \ref{2D Example Cooksmembrane - VEMVoronoi Mesh} and \ref{2D Example Cooksmembrane - VEMQ2S Mesh}, respectively. The contour plots of the vertical displacement evolution for different deformation states $\{t= 0~s, 0.0001~s, 0.0002~s, 0.00035~s, 0.00055~s, 0.00065~s \}$ are sketched in Figure \ref{CM-DM}. The nonlinear behavior is clearly observed in the deformation process due to the dynamic effects at finite strains.

\begin{figure}[!t]
	\begin{minipage}{0.5\textwidth}
		\hspace{30mm}
		\centering
		\begin{subfigure}[c]{0.5\textwidth}
			\scalebox{1}{\includegraphics{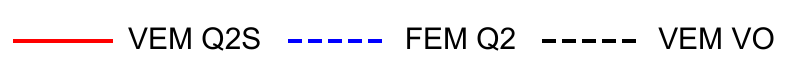}}
			\label{2D Example Cooksmembrane - Legend}
		\end{subfigure}
	\end{minipage}
\vspace{-2mm}

\begin{subfigure}[c]{0.5\textwidth}
		\scalebox{0.6}{\includegraphics{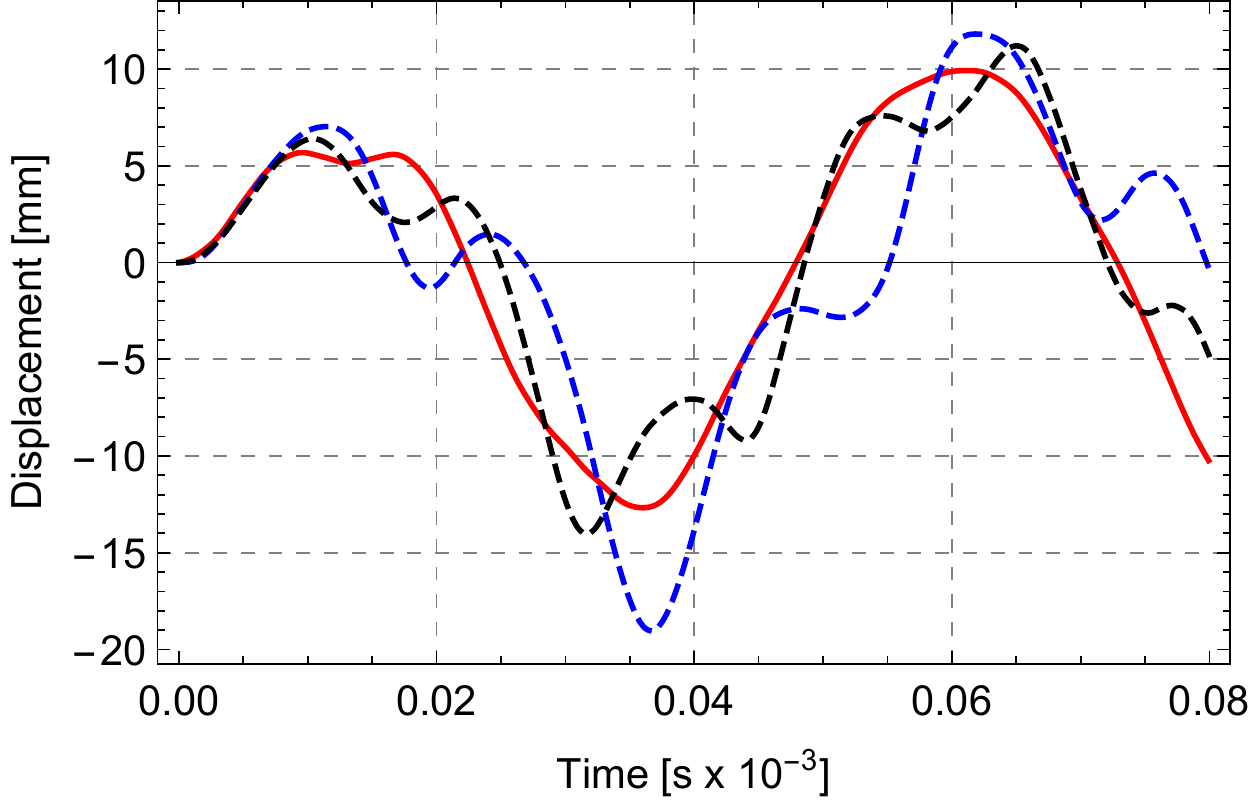}} 
		\vspace{-1mm}
		\subcaption{N=1}
		\label{2D Example Cooksmembrane - Comparison N1}
	\end{subfigure}
	\begin{subfigure}[c]{0.5\textwidth}
		\scalebox{0.6}{\includegraphics{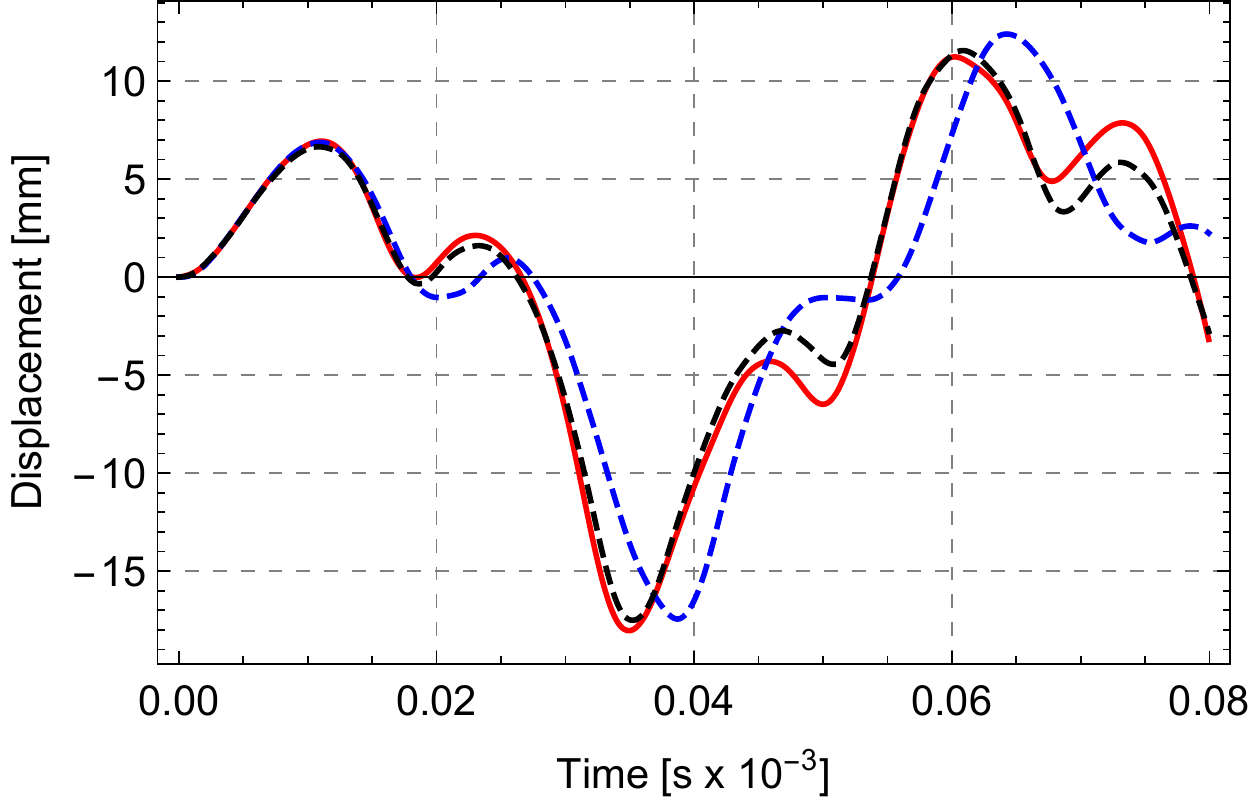}}
		\vspace{-1mm} 
		\subcaption{N=2}
		\label{2D Example Cooksmembrane - Comparison N2}
	\end{subfigure}
	\vspace{6mm}
	\begin{subfigure}[c]{0.5\textwidth}
		\scalebox{0.6}{\includegraphics{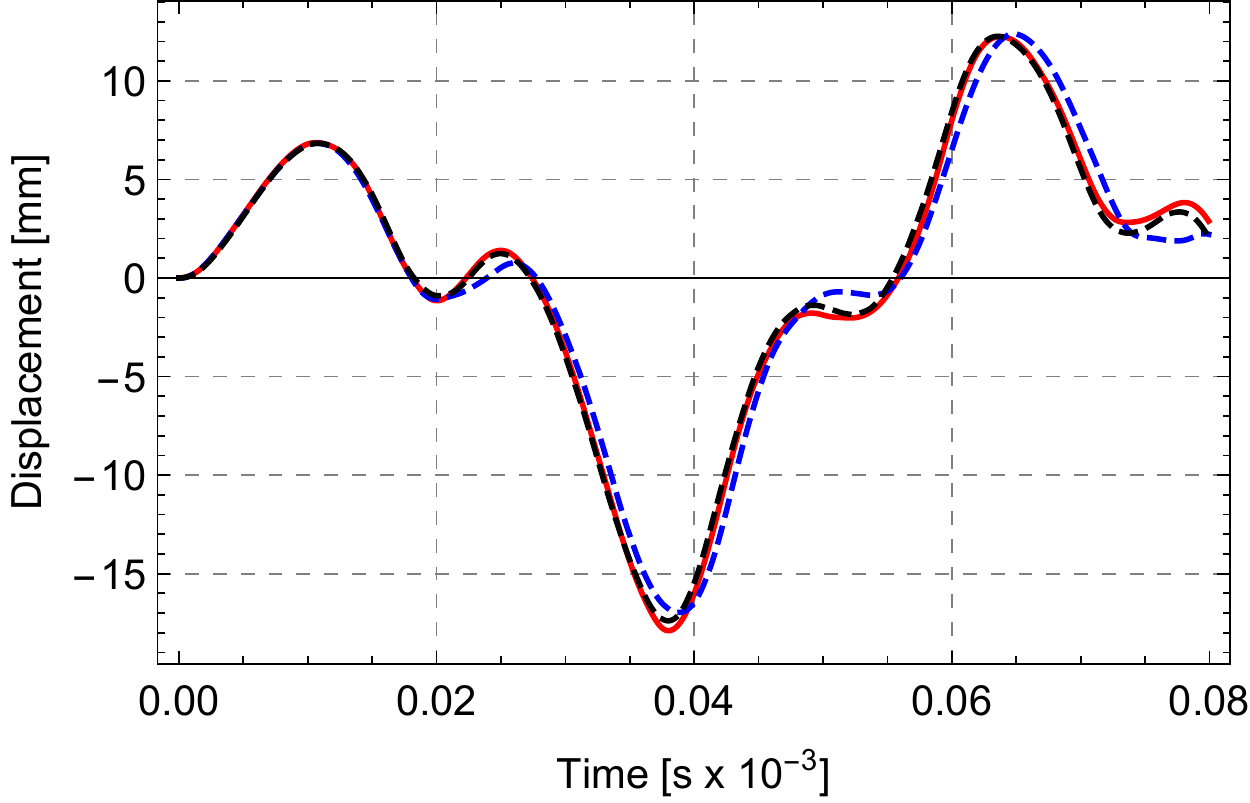}} 
		\vspace{-1mm}
		\subcaption{N=3}
		\label{2D Example Cooksmembrane - Comparison N3}
	\end{subfigure}
	\begin{subfigure}[c]{0.5\textwidth}
		\scalebox{0.6}{\includegraphics{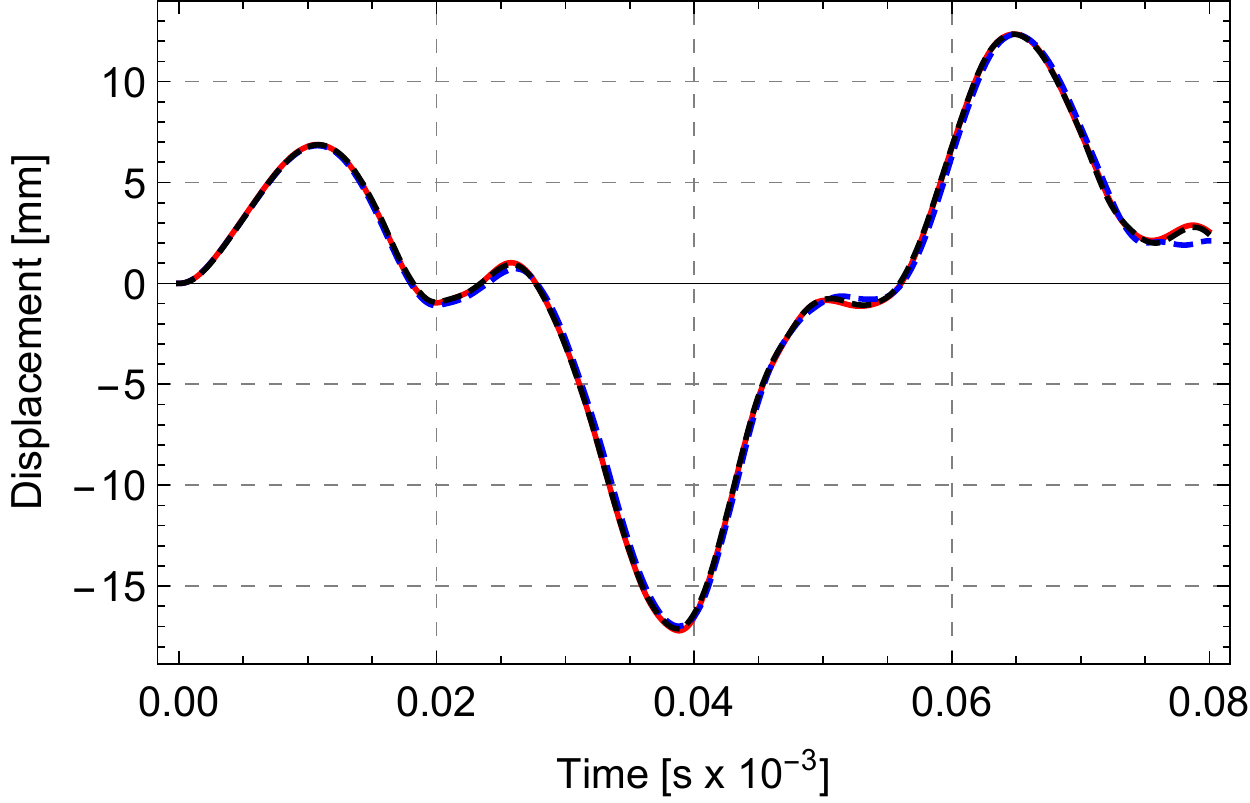}} 
		\vspace{-1mm}
		\subcaption{N=4}
		\label{2D Example Cooksmembrane - Comparison N4}
	\end{subfigure}
	\vspace{-2mm}
	\caption{ Convergence Study - Displacement over time response for 2D Example - Cook's membrane. Element division $2^N$, where $N$ increases from \textbf{(a)} to \textbf{(d)}.}
%
	\label{2D Example Cooksmembrane - Comparison}
\end{figure}

Figure \ref{2D Example Cooksmembrane - Comparison} shows a mesh refinement study with the element division of $2^N$ for N=1,2,3,4. For N=3 and higher the solution converges. {A comparison with FEM depicts that the results are in a very good agreement.} 

This study shows that again, that the evaluation of the integral in \eqref{eq:Integral_HTH} at the element centroid is absolutely sufficient to compute the mass-matrix.
\subsection{3D Boundary value problems}
\subsubsection{Wave propagation in a bar}
\begin{figure}[!t]
	\begin{minipage}{0.2\textwidth}
		\hspace{0mm}
		\centering
		\begin{subfigure}[!h]{0.5\textwidth}
			\scalebox{1}{\includegraphics{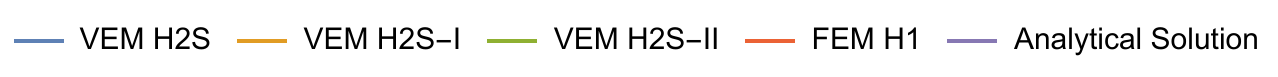}}
			\label{3D Example Longitudinal - Legend}
		\end{subfigure}
	\end{minipage}\\
	\begin{subfigure}[!h]{0.5\textwidth}
		\scalebox{0.62}{\includegraphics{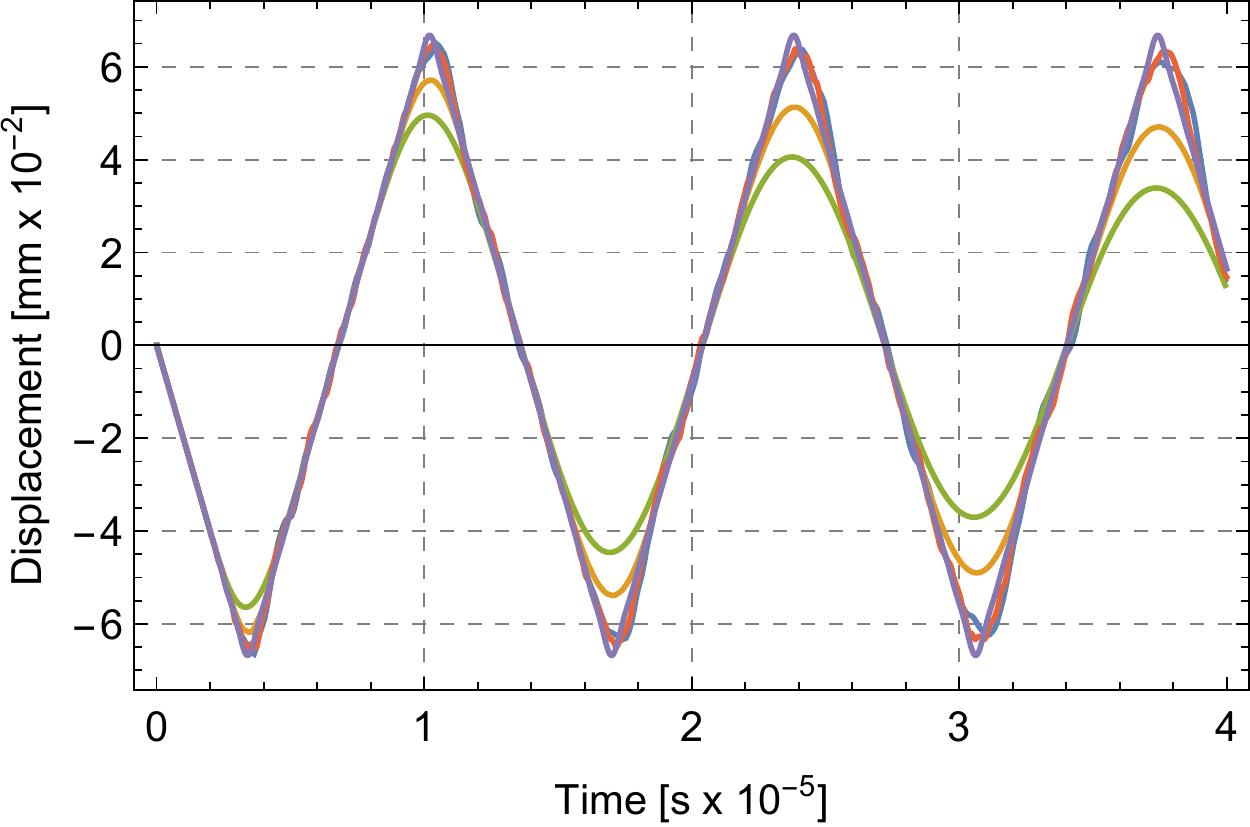}} 
		\subcaption{Response at $x=\ell$}
		\label{3D Example Longitudinal - Comparison 1}
	\end{subfigure}
	\begin{subfigure}[!h]{0.5\textwidth}
		\scalebox{0.62}{\includegraphics{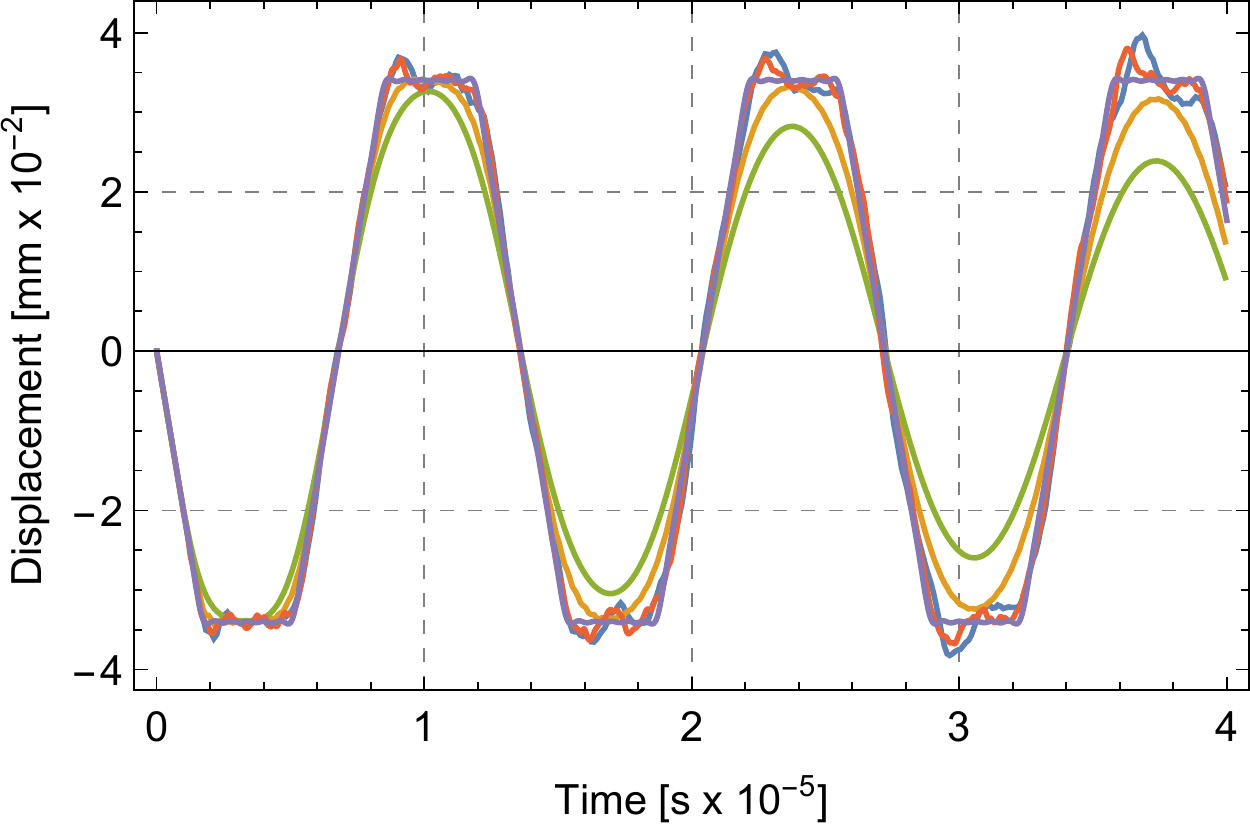}} 
		\subcaption{Response at $x=\ell/2$}
		\label{3D Example Longitudinal - Comparison 2}
	\end{subfigure}
	\caption{ Displacement over time response for 3D Example - Wave propagation in a bar that having an initial velocity of $v_0=20 m/s$.}
	\label{3D Example Longitudinal}
\end{figure}

The previously introduced 2D model of a bar is here extended to the third dimension. The length of the bar is set to $\ell=30mm$ and the height is equal to the width $h=b=5mm$. We apply an initial velocity of $v_0=20 m/s$ in longitudinal direction. The virtual element results are obtained using 400 elements, were the finite element results were obtained with 4320 elements. In this example we compare the virtual elements H2S, H2S-I and H2S-II with the finite element H1 and the analytical solution which was obtained for the 1D case in equation \eqref{eq:analyticalsolution1}. As already introduced before, the variable $\beta^{dyn}$ indicates how the mass-matrix is going to be evaluated. For $\beta^{dyn} =0$, the mass-matrix is calculated using only the projection part. Whereas for $\beta^{dyn} =1$ the computation of the mass-matrix is carried out using the stabilization part. Figure \ref{3D Example Longitudinal} depicts the displacement over time response at $x=\ell$ and $x=\ell/2$. The computation of the mass-matrix using VEM-H2S-I and VEM-H2S-II does not lead to accurate results. 
Whereas the computation using the projected part and evaluating the integral in \eqref{eq:Integral_HTH} at the element centroid (i.e. VEM-H2S) produces nearly the same results as the finite element H1 and the analytical solution.

\begin{figure}[b]
	\centering
	\scalebox{0.9}{\includegraphics{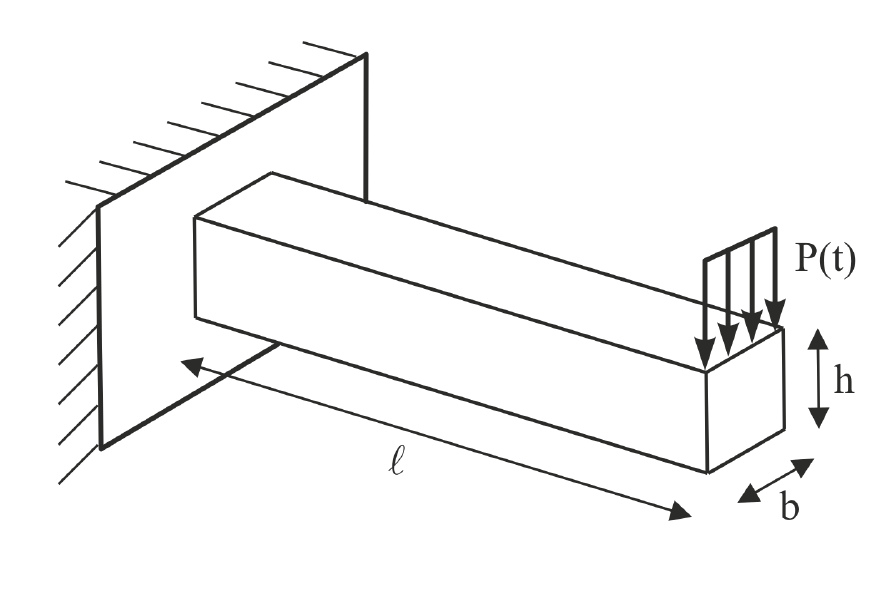}} 
	\caption{3D Example - Transversal vibration of a thick beam (Boundary value problem).}
	\label{3D Example Transversal - BVP}
\end{figure}

\subsubsection{Transversal vibration of a thick beam}

In this benchmark test a 3D cantilever beam is investigated. The geometric setup and the loading conditions of the specimen are depicted in Figure \ref{3D Example Transversal - BVP}. Here a line load is applied at the end of the beam with $P_{max} = 6  \ kN/mm$. The temporal course of the force is again given by a {\it half sine}, as shown in \ref{2D Example Transversal - Force}. Furthermore, similar to the 2D case, we set the beam length as $\ell=30 \ mm$ with equal height and width as $h=b=5 \ mm$. We compare the virtual element H1 with the finite elements H1 and H2. For this purpose a mesh refinement is employed from 8, 32, 128 to 1024 elements (N=1,2,3,4) for the finite element H1 and virtual element H1. The FEM H2 solution is computed with 3200 elements and can be seen as a {\it reference solution}.

\begin{figure}[!t]
	\begin{minipage}{0.2\textwidth}
		\hspace{0mm}
		\centering
		\begin{subfigure}[c]{0.5\textwidth}
			\centering
			\scalebox{1}{\includegraphics{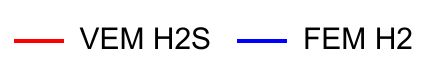}}
			\label{3D Example Transversal - Legend}
		\end{subfigure}
	\end{minipage}\\
	\begin{subfigure}[c]{0.5\textwidth}
		\scalebox{0.62}{\includegraphics{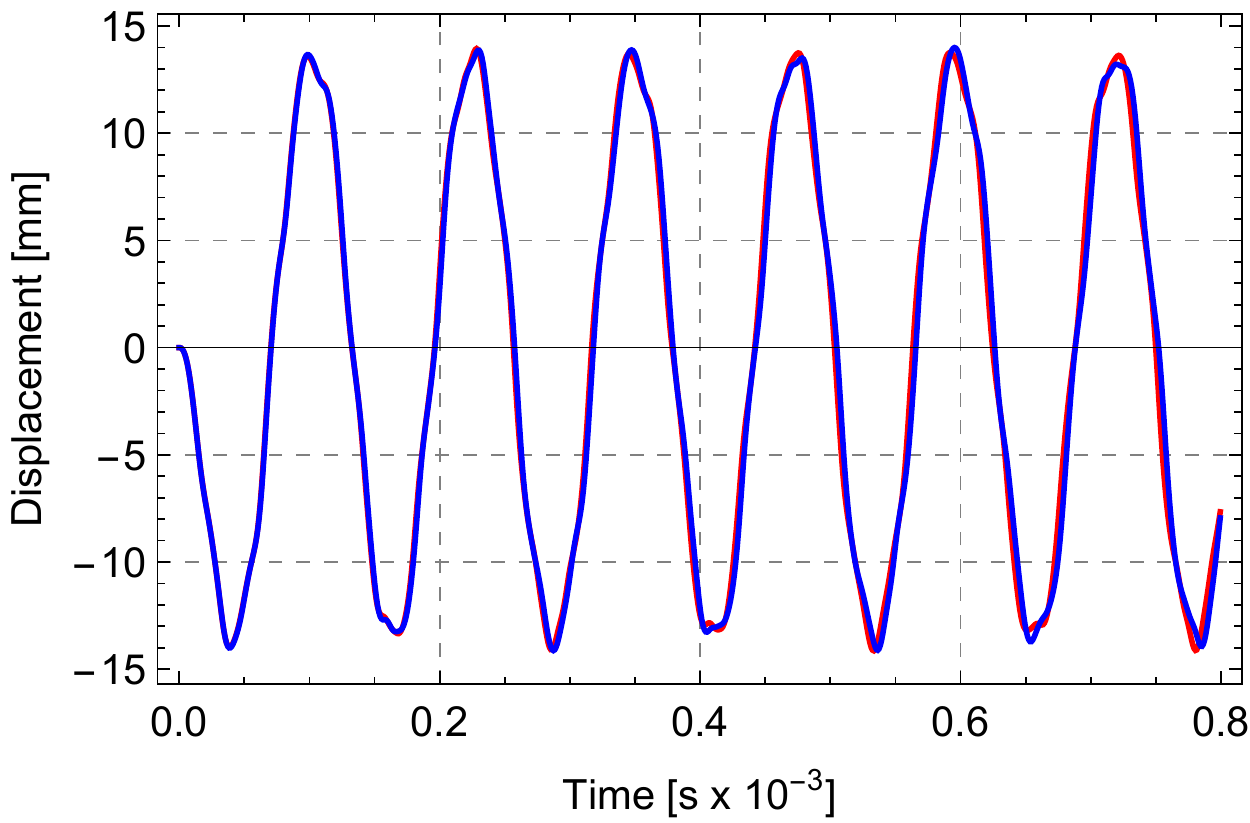}} 
		\subcaption{}
		\label{3D Example Transversal - Comparison 1}
	\end{subfigure}
	\begin{subfigure}[c]{0.5\textwidth}
		\scalebox{0.5}{\includegraphics{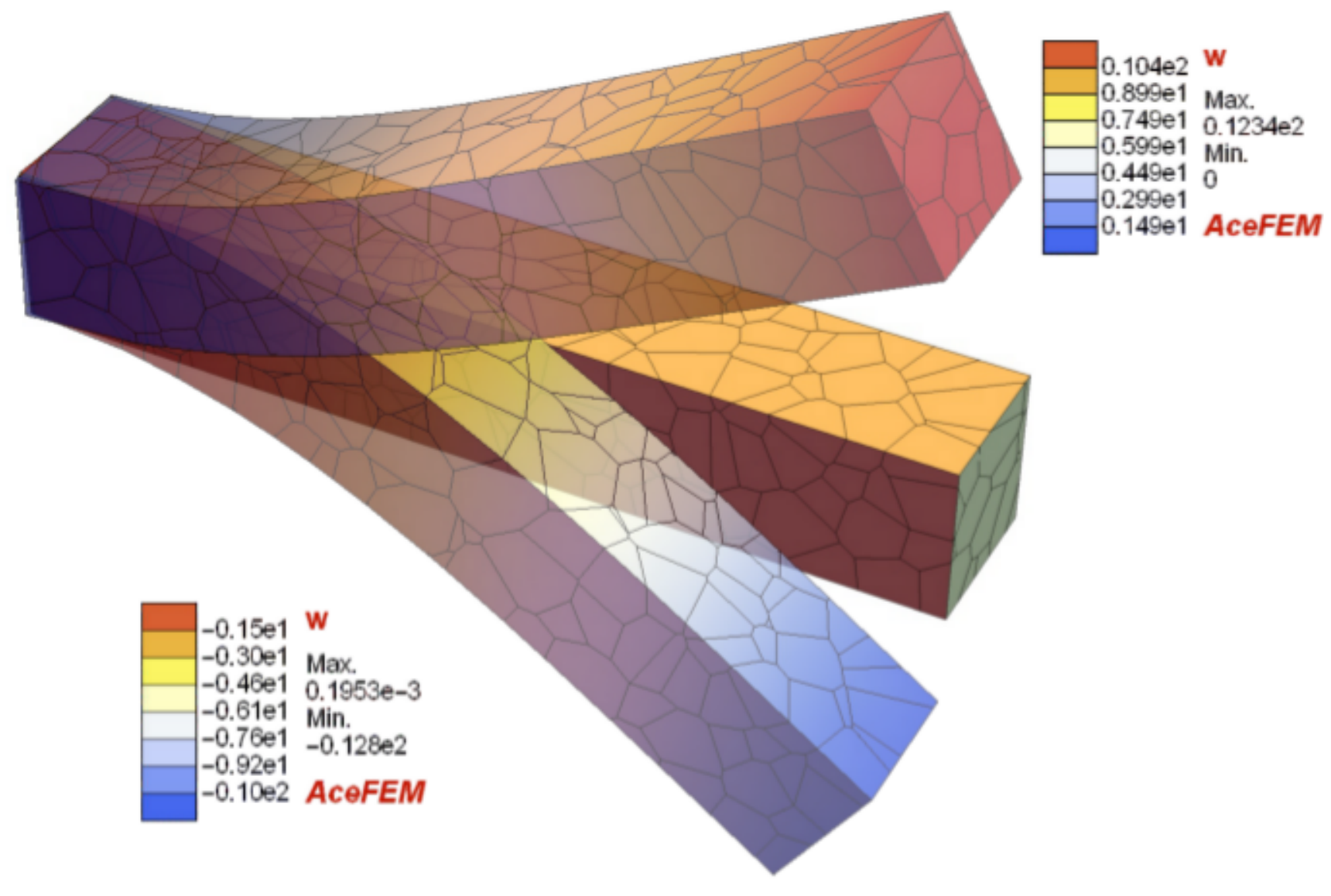}} 
		\subcaption{}
		\label{3D Example Transversal - Deformed Mesh}
	\end{subfigure}
	\caption{3D Example: Transversal vibration of a thick beam. (a) Displacement over time response at $x=\ell$ and (b) undeformed and maximal deformed mesh.}
\end{figure}

\begin{figure}[!h]
	\begin{minipage}{0.5\textwidth}
		\hspace{30mm}
		\centering
		\begin{subfigure}[c]{0.5\textwidth}
			\scalebox{1}{\includegraphics{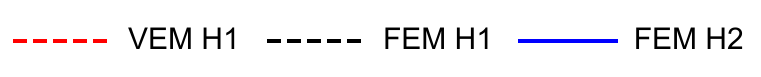}}
			\label{3D Example Transversal Beam Vibration - Legend}
		\end{subfigure}
	\end{minipage}
	
	\begin{subfigure}[c]{0.5\textwidth}
		\scalebox{0.62}{\includegraphics{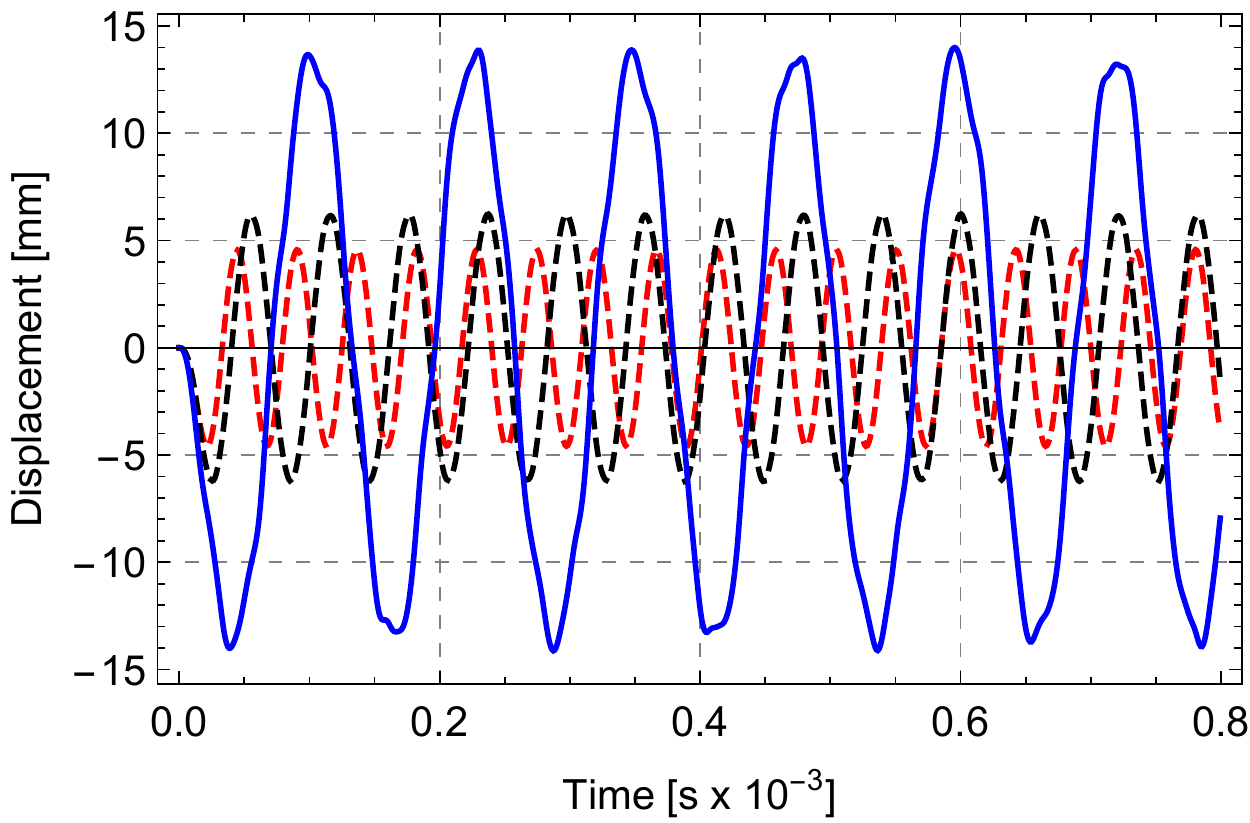}} 
		\vspace{-1mm}
		\subcaption{N=1}
		\label{3D Example Transversal Beam Vibration - Comparison N1}
	\end{subfigure}
	\begin{subfigure}[c]{0.5\textwidth}
		\scalebox{0.62}{\includegraphics{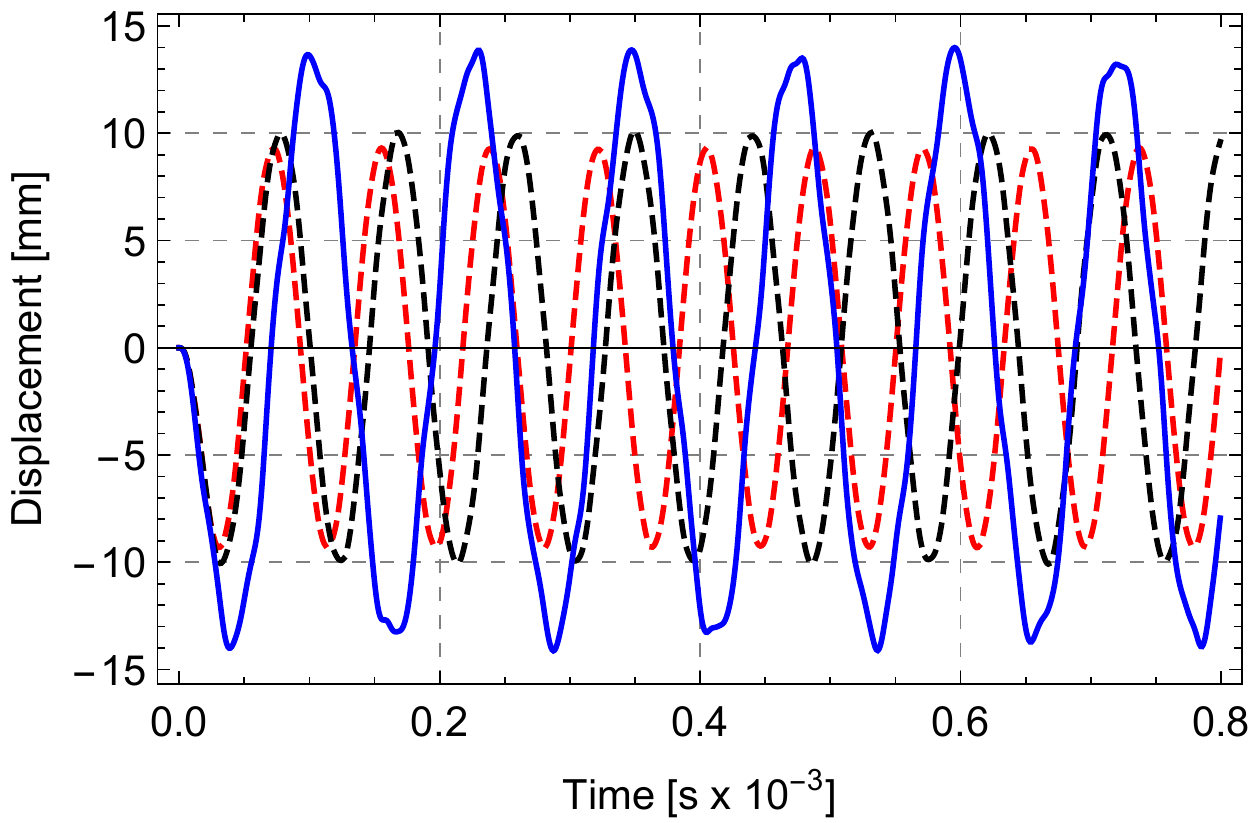}}
		\vspace{-1mm} 
		\subcaption{N=2}
		\label{3D Example Transversal Beam Vibration - Comparison N2}
	\end{subfigure}
	\vspace{6mm}
	
	\begin{subfigure}[c]{0.5\textwidth}
		\scalebox{0.62}{\includegraphics{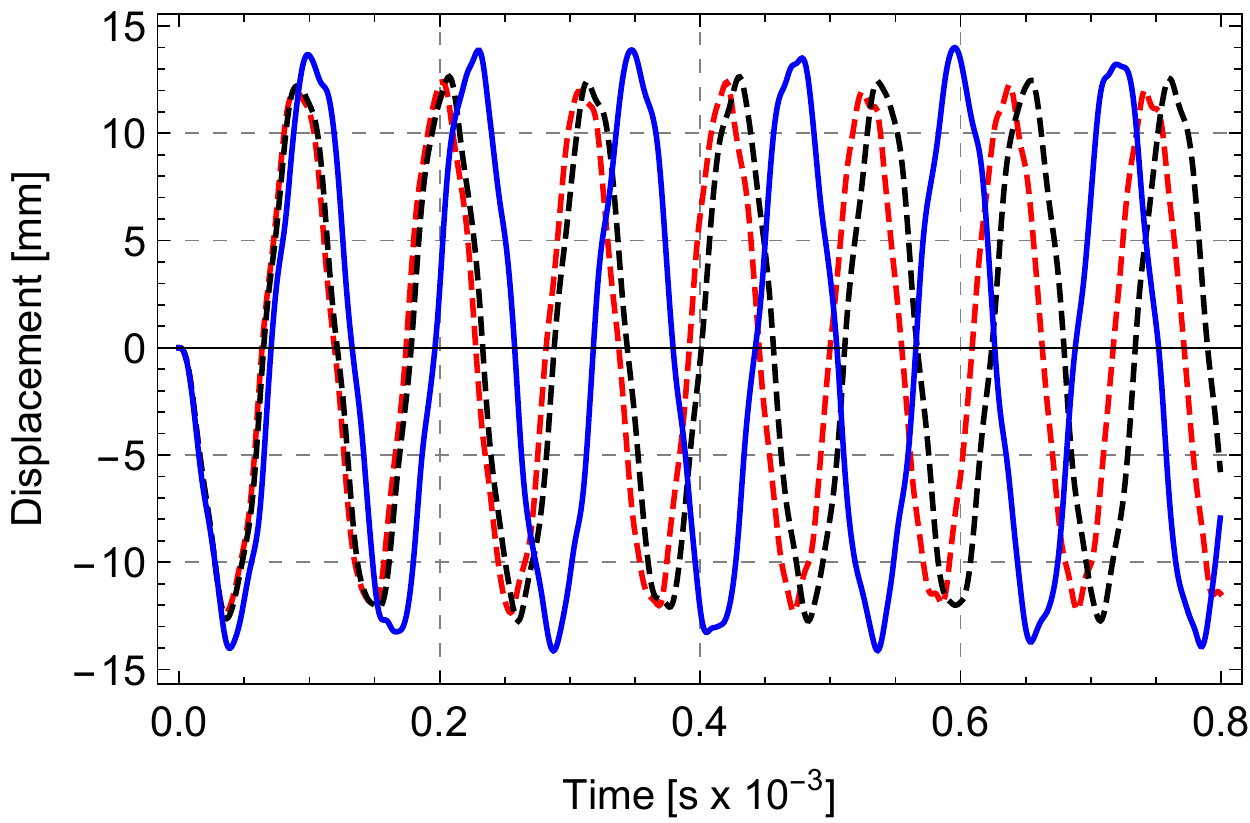}} 
		\vspace{-1mm}
		\subcaption{N=3}
		\label{3D Example Transversal Beam Vibration - Comparison N3}
	\end{subfigure}
	\begin{subfigure}[c]{0.5\textwidth}
		\scalebox{0.62}{\includegraphics{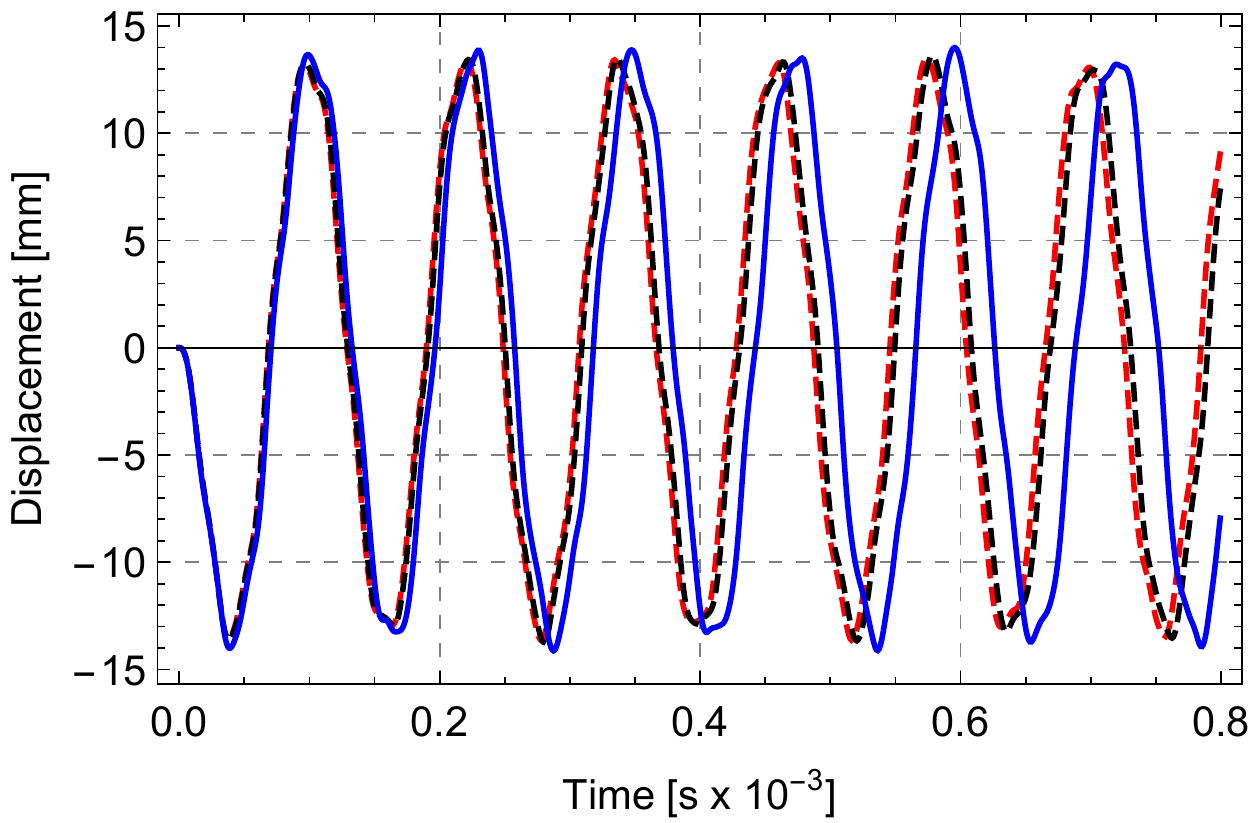}} 
		\vspace{-1mm}
		\subcaption{N=4}
		\label{3D Example Transversal Beam Vibration - Comparison N4}
	\end{subfigure}
	\caption{ Convergence Study - Displacement over time response for 3D Example - Transversal vibration of a thick beam. Element division $2^N$, where $N$ increases from \textbf{(a)} to \textbf{(d)}.}
	\label{Transversal Beam Vibration - Mesh Refinement}
\end{figure}

The maximum deformation state is sketched in Figure \ref{3D Example Transversal - Deformed Mesh}, representing the deflection $w$. Here the nonlinear behavior is clearly observed due to the dynamic effects at finite strains. Figure \ref{Transversal Beam Vibration - Mesh Refinement} illustrates the displacement over time response at $x=\ell$ for the mesh refinement study. This response is plotted for the center of the cross section. We observe that both, VEM and FEM results are converging to the reference solution for increasing number of elements. Still there is a shift with increasing time, this is due to the less accuracy of VEM/FEM H1 element compared with the FEM H2 quadratic ansatz function.

Additionaly, we employ the virtual element VEM H2S computed with 256 elements in Figure \ref{3D Example Transversal - Comparison 1} and compare it with the reference solution (FEM H2 with 3200 elements). It is interesting to note that despite the use of linear ansatz functions VEM H2S produces nearly the same results as the reference solution. This is due to the fact that the stabilization uses the bending modes. In conclusion, the presented formulation depicts very good results also in the 3D case.

This test also confirms that the evaluation of the integral in \eqref{eq:Integral_HTH} for the computation of the mass-matrix only at the centroid of the polygon/polyhedra is absolutely enough to get satisfying results.

\subsubsection{Vibration of a thick plate}
\begin{figure}[!t]
	\centering
	\scalebox{0.7}{\includegraphics{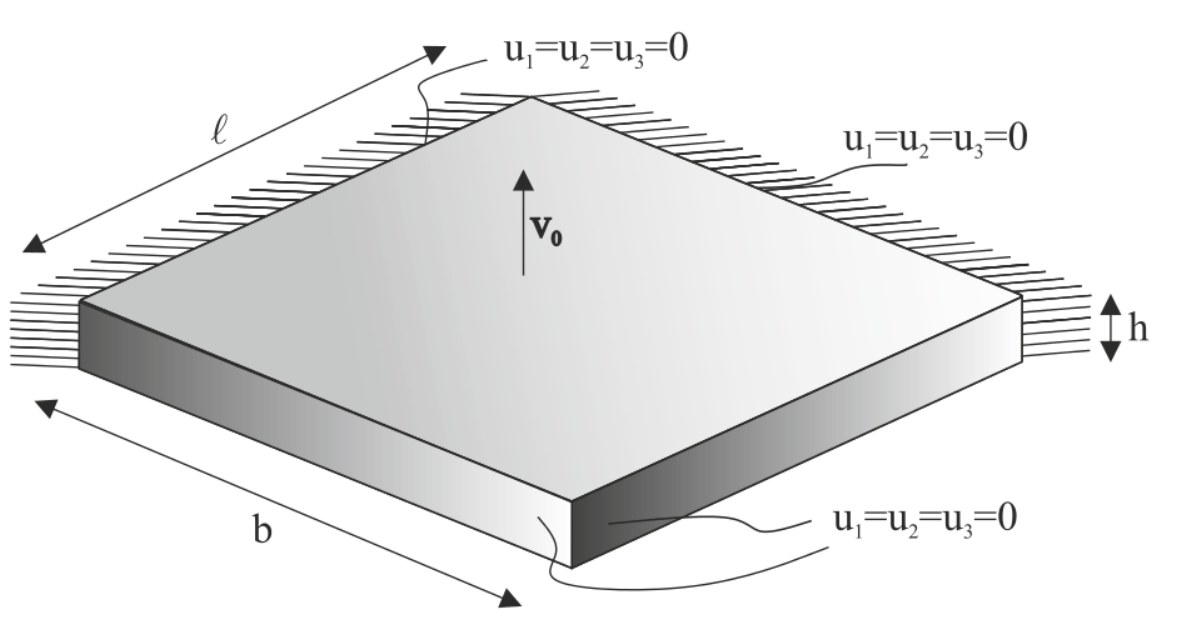}} 
	\caption{Vibration of a thick plate (Boundary value problem)}
	\label{3D Example Plate - BVP}
\end{figure}

The last example is related to the vibration of a thick plate which is discretized using three-dimensional elements. The plate has a length $\ell=30 \ mm$, a thickness $h=5 \ mm$ and a width $b=30 \ mm$ as shown in Figure \ref{3D Example Plate - BVP}. We further set an initial condition such that the initial velocity of the plate is set to {$v_{0}=200 m/s$} in the $z$-direction, see Figure \ref{3D Example Plate - BVP}.

\begin{figure}[!b]
	\begin{subfigure}[t]{0.32\textwidth}
		\scalebox{0.57}{\includegraphics{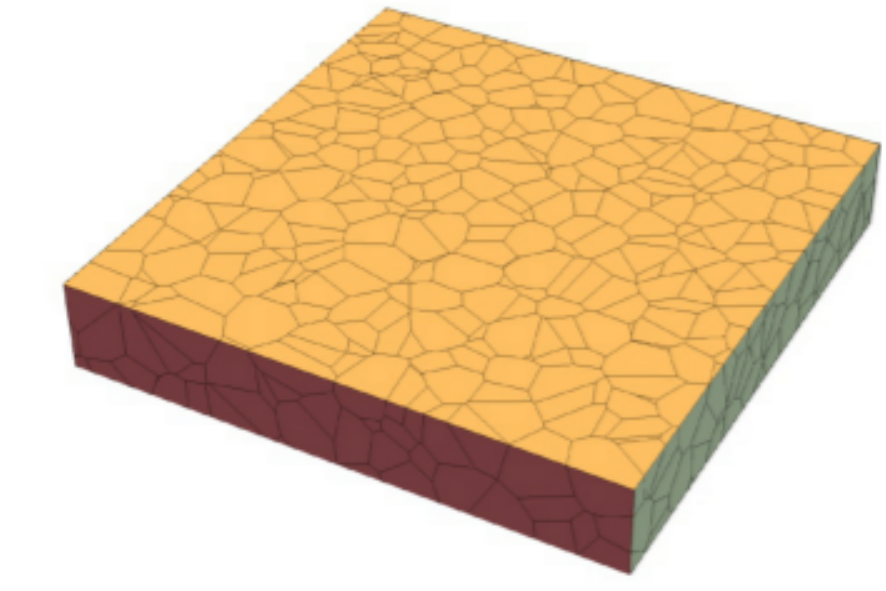}}
		\label{svar009 Deformed Mesh 001}
		\subcaption{t=0s}
	\end{subfigure}
	\begin{subfigure}[t]{0.32\textwidth}
		\scalebox{0.57}{\includegraphics{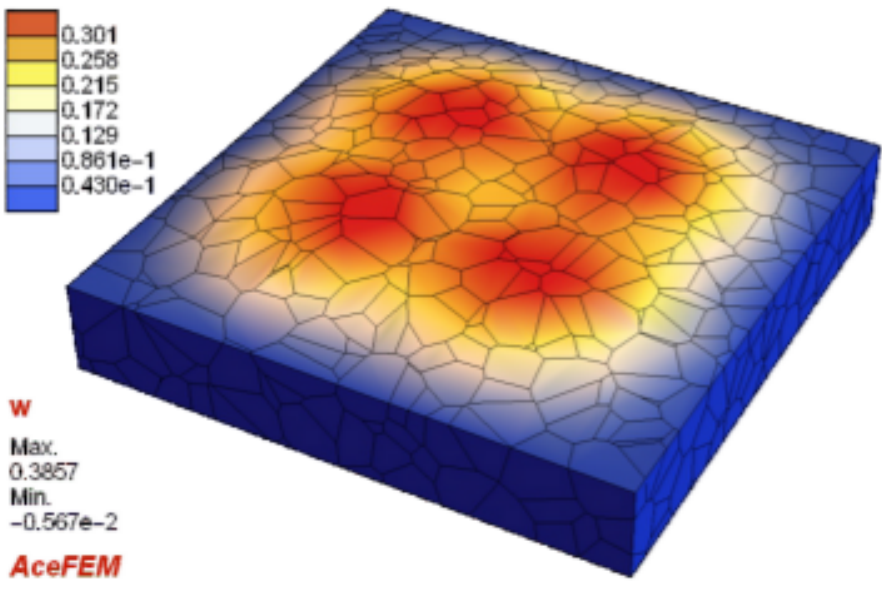}} 
		\label{svar009 Deformed Mesh 002}
		\subcaption{t=0.0000014798}
	\end{subfigure}
	\begin{subfigure}[t]{0.32\textwidth}
		\scalebox{0.57}{\includegraphics{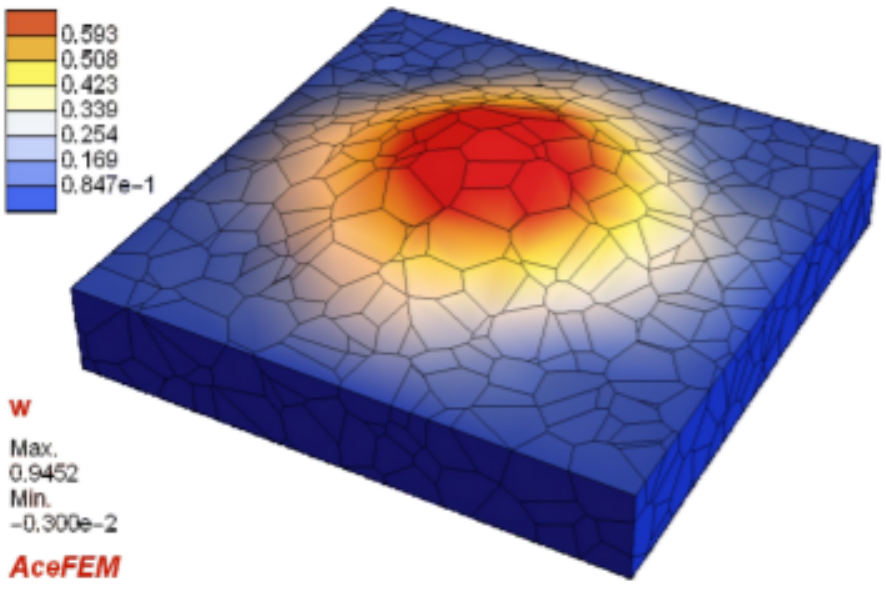}} 
		\label{svar009 Deformed Mesh 003}
		\subcaption{t=0.00000290479}
	\end{subfigure}\\[5mm]
	\begin{subfigure}[t]{0.32\textwidth}
		\scalebox{0.57}{\includegraphics{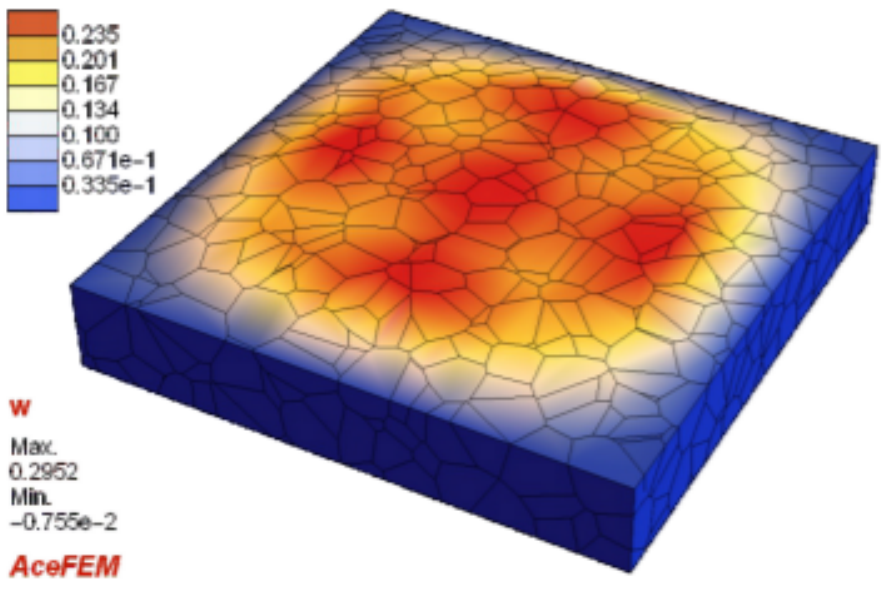}} 
		\label{svar009 Deformed Mesh 004}
		\subcaption{t=0.00000535609}
	\end{subfigure}
	\begin{subfigure}[t]{0.32\textwidth}
		\scalebox{0.57}{\includegraphics{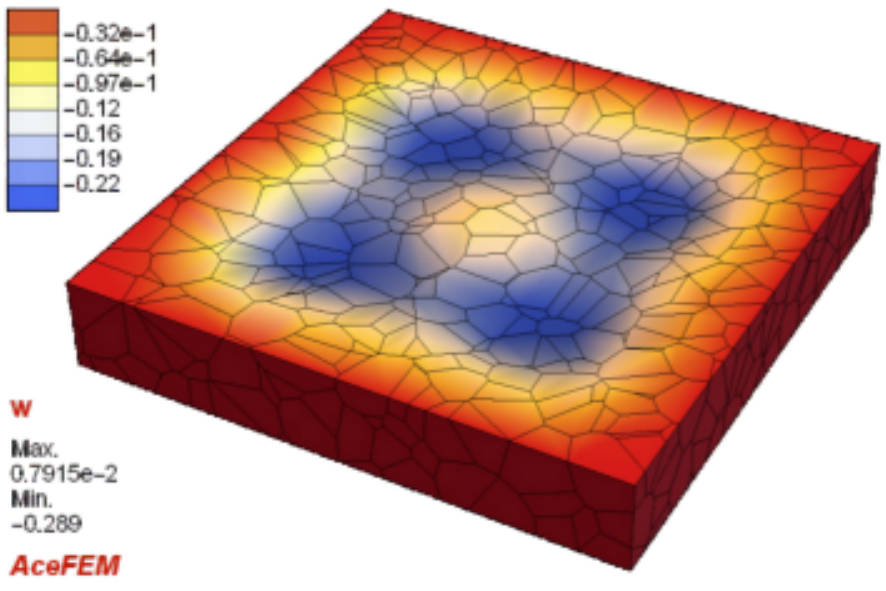}} 
		\label{svar009 Deformed Mesh 005}
		\subcaption{t=0.000007597}
	\end{subfigure}
	\begin{subfigure}[t]{0.32\textwidth}
		\scalebox{0.57}{\includegraphics{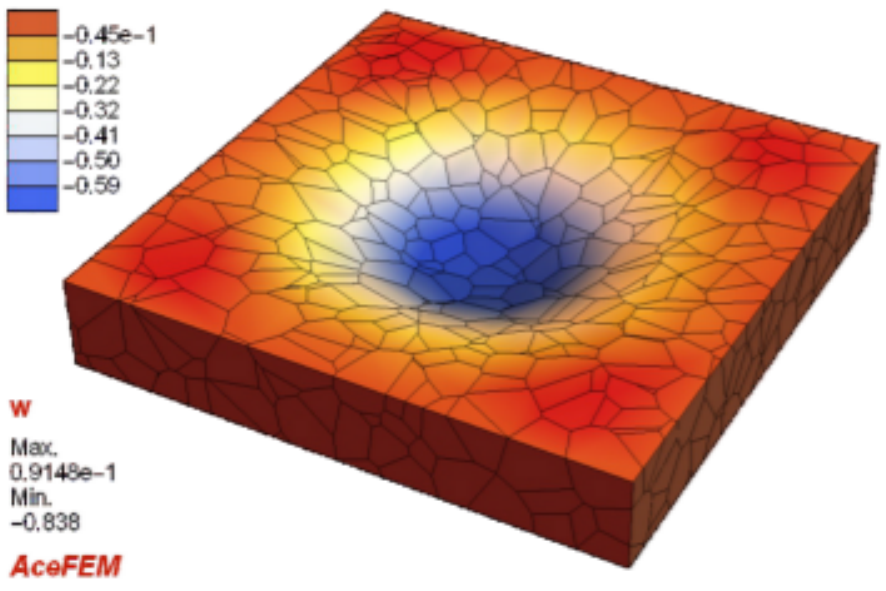}} 
		\label{svar009 Deformed Mesh 006}
		\subcaption{t=0.00000989465}
	\end{subfigure}\\[4mm]
	\caption{3D Example Vibration of a thick plate - {Evolution of the displacement in the $z$-direction for different deformation states.}}
		\label{plate-defos}	
\end{figure}

\begin{figure}[!t]
\centering
\includegraphics{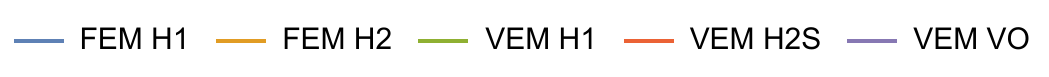}\\[2mm]
\includegraphics[width=0.5\textwidth]{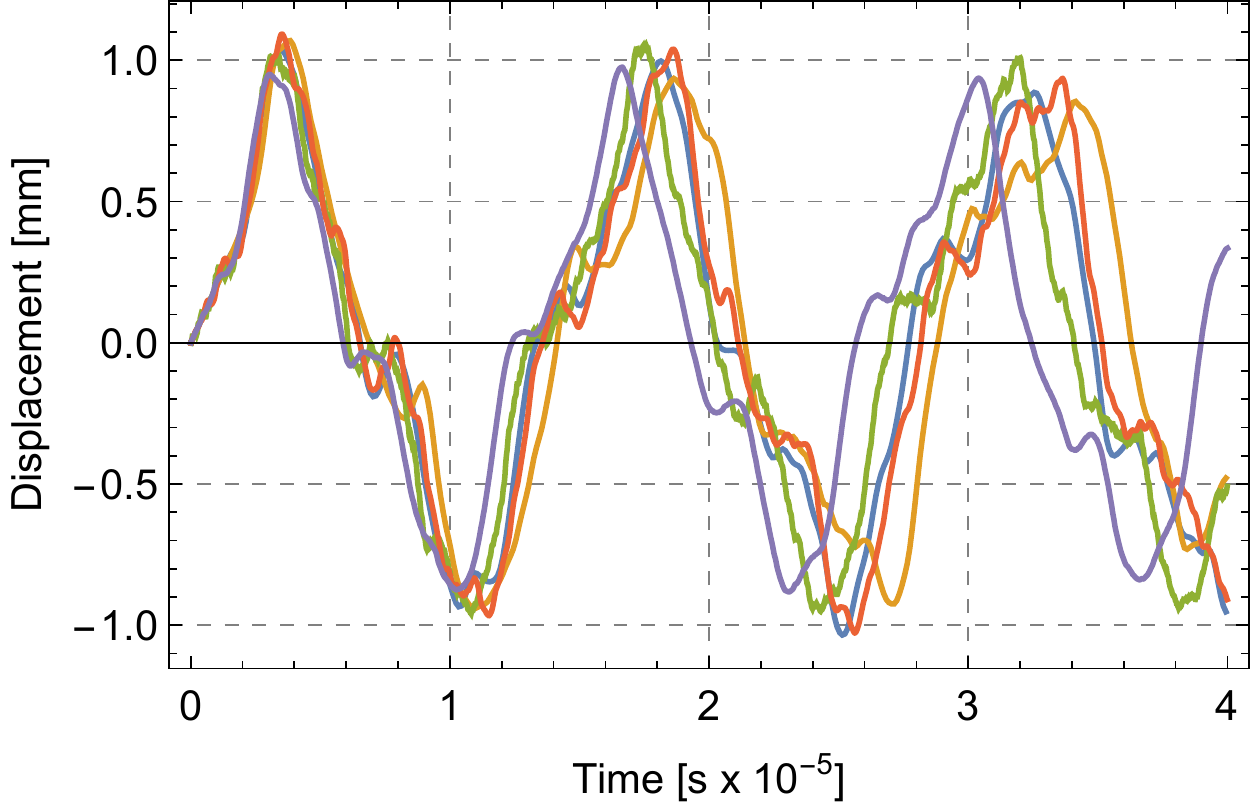}\\[2mm]
\caption{3D example - Vibration of a thick plate: Vertical displacement over time response at the center of the plate.}
\label{3DExamplePlate-Comparison}
\end{figure}

Figure \ref{plate-defos} demonstrates the evolution of the displacement in the $z$-direction for different deformation states using the VEM-VO element. Herein a nonlinear response undergoing large deformation is observed due to the elastodynamic behavior. In \ref{3DExamplePlate-Comparison} the vertical displacement over time response is plotted at the center of the plate at the thickness ${z=h/2}$. We observe a good match between the virtual element results and the finite element results. Here we used in addition to regular shaped elements, Voronoi elements which have an arbitrary number of nodes and element shapes. The computation is performed with 1024 virtual elements of type H1, H2S and VO and the finite element H1. We used 6400 elements for FEM-H2 as a reference solution. Again one can observe that the computation of the mass-matrix using only the projection part and evaluating the integral in \eqref{eq:Integral_HTH} at the centroid of the element yields accurate results. 


\section{Summary and Conclusions}
\label{Summary}

An efficient low order virtual element formulation was developed for nonlinear elastodynamics within this work. The virtual element results show a very good match with finite element results and analytical results. Arbitrary shaped elements with a various number of nodes could be used successfully for the simulations.\\
The computation of the mass-matrix was performed using the projection part and does not need any stabilization. This is valid only for problems, where the governing equations have no reaction terms \cite{beirao2013c, Ahma13}, which was the case in this work. To compute the integral in equation \eqref{eq:Integral_HTH}, the argument can be evaluated at the element centroid. This is sufficient accurate as shown in the examples. Hence, there is no need to perform any sub-triangulation of the element or use the moment of areas in \eqref{eq:moments of area} for the computation of the mass-matrix.  
\\
The extension of VEM to other applications open a wide range of new research directions such as dynamic elasto-plasticity, contact or impact.

\bigskip

{\bf Acknowledgments.}
The authors gratefully acknowledges support for this research by the ``German Research Foundation'' (DFG) in (i) the {\sc collaborative research center} {\bf CRC 1153}, (ii) the {\sc Priority Program {\bf SPP 2020}} and (iii) {\sc Priority Program {\bf SPP 1748}}.
\bibliography{nlbuch,PW,parallel,fem02,Lit_VEM_cont,BH_bibliography,MC_bibliography,thermo}
\bibliographystyle{ieeetr}

\end{document}